\documentclass[3p,10pt]{elsarticle}
\usepackage{amsfonts,amssymb,amsmath,amsthm,indentfirst}
\usepackage{bm}
\usepackage{pgfplotstable,booktabs}
\usepackage{graphicx}
\usepackage{caption}
\newtheorem{theorem}{Theorem}
\newtheorem{lemma}{Lemma}
\newtheorem{remark}{Remark}

\def\ve{\varepsilon}

\def\vp{\varphi}
\def\leq{\leqslant}
\def\geq{\geqslant}
\def\onetwelve{{\textstyle\frac{1}{12}}}
\def\half{{\textstyle\frac12}}
\def\*#1{\mathbf{#1}}
\def\urho{\underline{\rho}}

\begin{document}
\begin{frontmatter}
\title{On 
%properties of 
a semi-explicit 
%in time 
fourth-order vector compact scheme for the multidimensional acoustic wave equation}
\author[az]{Alexander Zlotnik\corref{cor1}}
\address[az]{Department of Mathematics, Higher School of Economics University, Pokrovskii bd. 11, 109028 Moscow, Russia}
\cortext[cor1]{Corresponding author}
\ead{azlotnik@hse.ru}
\author[az]{Timofey Lomonosov}
\ead{tlomonosov@hse.ru}

\begin{abstract}
\noindent We deal with an initial-boundary value problem for the multidimensional acoustic wave equation, with the variable speed of sound.
For a three-level semi-explicit in time higher-order vector compact scheme, we prove stability and derive 4th order error bound in the enlarged energy norm.
This scheme is three-point in each spatial direction, and it exploits additional sought functions which approximate 2nd order non-mixed spatial derivatives of the solution to the equation.
At the first time level, a similar two-level in time scheme is applied, with no derivatives of the data.
No iterations are required to implement the scheme.
We also present results of various 3D numerical experiments 
that demonstrate a very high accuracy of the scheme for smooth data, its advantages in the error behavior over the classical explicit 2nd order scheme for nonsmooth data as well and an example of the wave in a layered medium initiated by the Ricker-type wavelet source function.
\end{abstract}

\begin{keyword}
acoustic wave equation, semi-explicit vector scheme, higher-order compact scheme, conditional stability, error bound
\end{keyword}
\end{frontmatter}
\par AMS Subject Classification: 65M06, 65M12, 65M15.

\section{Introduction}
\label{s:intro}

This paper continues a study of the semi-explicit in time fourth-order vector compact scheme for the multidimensional acoustic wave equation, with the variable speed of sound, that has recently been constructed and investigated in \cite{ZArxiv21,ZL24}.
A nonstandard element of the scheme is the use of additional sought functions that approximate the 2nd order nonmixed spatial derivatives of the solution; that element leads to a simple and direct implementation of the scheme.
Its conditional stability in the strong energy norm and the corresponding error bound of the order 3.5 have been proved in \cite{ZL24}. 
%In a simpler case of 
For the wave equation (having
%with 
constant coefficients), stability in the standard energy norm and the corresponding 4th order error bound have been proved as well, but the argument is not extended to the case of variable coefficients where the 4th order bound has not been proved.

\par In this paper, we apply another approach based on a theorem of stability in an enlarged weak energy norm for an abstract three-level method together with symmetrization of the scheme being studied. 
This makes it possible to prove the expected stability in the enlarged (standard) energy norm and the corresponding 4th order error bound. 
Also, in addition to 2D numerical experiments in \cite{ZL24}, we accomplish new 3D ones.  
First, we show very small errors and their 4th order behavior for smooth travelling wave-type solutions.  
Second, we demonstrate essential advantages in error and its convergence rates in the mesh $L^2$, $H^1$ and energy norms of the present scheme over the classical explicit 2nd order scheme in the case of nonsmooth data as well. 
Such important effect is observed namely in the hyperbolic case, and not in the elliptic or parabolic ones. 
Note that, in 1D case, the practical convergence rates in the case of nonsmooth data were studied for some 2nd and 4th order schemes in \cite{ZK18,ZK20}. 
Finally, we give an example of a wave propagation in a three-layer medium, with a piecewise constant speed of sound, initiated by a smoothed Ricker-type wavelet source. 

\par The semi-explicit vector compact scheme was originally suggested for the 2D wave equation in \cite{JG20}.
Its generalization to the $n$-dimensional wave and acoustic wave equations, $n\geq 1$, has been given in \cite{ZArxiv21} (including the case of nonuniform meshes in space and time).
Recall that the case $n\geq 4$ is of interest in theoretical physics, for example, see \cite{F17}.
The cases of the 3D wave and acoustic wave equations have been also considered in \cite{JG23}, where another spatial discretization in the latter case has been suggested; the case of mildly nonlinear wave equation has been covered as well. 
Notice that a rigorous justification of the last mentioned discretization is absent and seems to be a technically complicated task.

\par We should notice that other compact 4th-order schemes were developed much earlier, and a vast literature is devoted to them.
The most standard of them were derived based mainly on an analysis and approximation of the leading term of the approximation error of second order schemes, in particular, see \cite{BTT18,MA22,P98,STT19,ZK20,ZC23}, etc. 
For the acoustic wave equation, rigorous stability theorems and corresponding 4th order error bounds for schemes of such type have recently been proved in \cite{ZC23}. 
But such schemes are implicit, and their effective implementation for 
%in the case of 
the acoustic wave equation needs application of iterative methods that entails increasing of computational costs.

\par A natural development of such compact schemes led to their various alternating-direction implicit (ADI) versions, see such 4th order schemes, in particular, in \cite{CL75,DZ13,DZ13b,LS13,Liao14,ZZLC10,ZK20}
and, in the case of the acoustic wave equation, in \cite{LLL19,LYDH18,ZC22}.  
To implement them, one needs to solve only 1D systems of linear algebraic equations with tridiagonal matrices; thus, no iterative methods are required, that is essentially simpler.
%Their implementation is much simpler, direct and reduces to solving only collections of 1D systems of linear algebraic equations with tridiagonal matrices.
But, for the acoustic wave equation, rigorous stability theorems and corresponding 4th order error bounds face
essential difficulties and are not available up to the moment.
%technical difficulties and are still absent.
Note that implicitness and conditional stability are immanent to all three types of compact 4th order schemes listed above.
%Note that all listed three types of 4th order compact schemes are implicit and conditionally stable.

\par We also mention that, for the wave-type equations, there exist a lot of quite different higher-order methods, in particular, see \cite{BB79,BDE21,CDI20,CFHJSS18,C02,HLZ19,SKWK18} and references therein.

\par Content of the paper is as follows.
In Section \ref{expl vector scheme}, we state the initial-boundary value problem (IBVP) for the acoustic wave equation, with nonhomogeneous Dirichlet boundary conditions, briefly recall how to construct the semi-explicit vector compact scheme to solve it 
%the derivation of the semi-explicit vector compact scheme for it 
and discuss its fast implementation.
For this scheme, in the next Section \ref{stab and error bounds}, %theorems on 
a stability theorem in the enlarged energy norm and
the corresponding theorem on 
the 4th order error bound are derived.
%in this norm are proved.
In the last Section \ref{numer results}, we give results of various 3D numerical experiments demonstrating the very high accuracy of the scheme for smooth data, its advantages in the error  behavior over the classical explicit 2nd order scheme for nonsmooth data as well and an example of the wave in a layered medium initiated by the Ricker-type wavelet source function.
\section{The multidimensional acoustic wave equation and a semi-explicit  4th-order vector compact scheme to solve it}
%\section{The $n$-dimensional acoustic wave equation and a semi-explicit in time 4th-order vector compact scheme for it}
\label{expl vector scheme}

\par 
In this paper, we treat
%The following 
an IBVP (initial-boundary value problem) for the $n$-dimensional acoustic wave equation %under the nonhomogeneous Dirichlet boundary condition 
%is treated
\begin{gather}
\rho(x)\partial_t^2u(x,t)-Lu(x,t)=f(x,t),
\ \ L:=a_1^2\partial_1^2+\ldots+a_n^2\partial_n^2,\ \ \text{in}\ \ Q_T=\Omega\times (0,T);
\label{hyperb2eq}\\
 u|_{\Gamma_T}=g(x,t);\ \ u|_{t=0}=u_0(x),\ \ \partial_tu|_{t=0}=u_1(x),\ \ x\in\Omega=(0,X_1)\times\ldots\times(0,X_n),
\label{hyperb2ibc}
\end{gather}
where the nonhomogeneous Dirichlet boundary condition is used.
Here $0<\urho\leq\rho(x)\leq\bar{\rho}$ and $a_1>0,\ldots,a_n>0$ are constants (they are taken different like in \cite{ZK20,ZC23,ZL24}), $x=(x_1,\ldots,x_n)$, $n\geq 1$.
In addition, $\partial\Omega$ and $\Gamma_T=\partial\Omega\times (0,T)$ are the boundary of $\Omega$ and the lateral surface of $Q_T$, respectively.
%Also $\partial\Omega$ is the boundary of $\Omega$ and $\Gamma_T=\partial\Omega\times (0,T)$ is the lateral surface of $Q_T$.
Let $a_{\max}:=\max_{1\leq k\leq n}a_k$.
We first assume that $\rho$ and $u$ are smooth (the nonsmooth examples will be taken in Section \ref{numer results}).

\par We introduce the uniform mesh $\overline\omega_{h_t}$ with the nodes $t_m=mh_t$, $0\leq m\leq M$ and the step $h_t=T/M>0$ on 
%a segment 
$[0,T]$; here $M\geq 2$.
Let $\omega_{h_t}=\{t_m\}_{m=1}^{M-1}$.
%be its internal part.
We define the mesh averages including the Numerov one, difference operators and the %operator of 
summation operator in time levels
%in $t$
with the variable upper limit 
%in $t$
\begin{gather*}
 \bar{s}_ty=\half(\check{y}+y),\ \
 s_{tN}y=\tfrac{1}{12}(\check{y}+10y+\hat{y}),\ \
 \bar{\delta}_ty=\frac{y-\check{y}}{h_t},\ \
 \delta_ty=\frac{\hat{y}-y}{h_t},\ \
\Lambda_ty=\delta_t\bar{\delta}_ty=\frac{\hat{y}-2y+\check{y}}{h_t^2},
\\
 I_{h_t}^my=h_t\sum_{l=1}^m y^l\ \ \text{for}\ \ 1\leq m\leq M,\ \ I_{h_t}^0y=0.
\end{gather*}
where $y^m=y(t_m)$, $\check{y}^{m}=y^{m-1}$ and $\hat{y}^{m}=y^{m+1}$.

\par We introduce the uniform mesh $\bar{\omega}_{hk}$  with the nodes $x_{ki}=ih_k$, $0\leq i\leq N_k$, with the step $h_k=X_k/N_k$ in $x_k$.
Denote by $\omega_{hk}=\{x_{ki}\}_{i=1}^{N_k-1}$ its internal part.
We define the difference operator and the Numerov average in $x_k$:
\begin{gather}
 (\Lambda_kw)_i:=\tfrac{1}{h_k^2}(w_{i+1}-2w_i+w_{i-1}),\ \
 s_{kN}w_i:=\tfrac{1}{12}(w_{i-1}+10w_i+w_{i+1})=(I+\onetwelve h_k^2\Lambda_k)w_i
\label{basic oper2}
\end{gather}
on $\omega_{h_k}$, where $w_i=w(x_{ki})$.

\par Further, we introduce the rectangular mesh $\bar{\omega}_h=\bar{\omega}_{h1}\times\ldots\times\bar{\omega}_{hk}$ in $\bar{\Omega}$,
where $h=(h_1,\ldots,h_n)$.
Let
%We denote by 
$\omega_h={\omega}_{h1}\times\ldots\times{\omega}_{hk}$ and $\partial\omega_h=\bar{\omega}_h\backslash\omega_h$.
%the internal part and boundary of $\bar{\omega}_h$.
We 
%also 
introduce the meshes $\omega_{\*h}:=\omega_h\times\omega_{h_t}$ in $Q_T$ and $\partial\omega_{\*h}=\partial\omega_h\times\{t_m\}_{m=1}^M$ on $\bar{\Gamma}_T$, with
%where 
$\*h=(h,h_t)$.

\par We denote by $H_h$ a Euclidean space of functions, with the inner product $(v,w)_h$ and norm $\|w\|_h=(w,w)_h^{1/2}$. For a linear operator $C_h$  acting in $H_h$, we denote by $\|C_h\|$ its norm.
Any such $C_h=C_h^*>0$ generates the norm $\|w\|_{C_h}:=(C_hw,w)_h^{1/2}=\|C_h^{1/2}w\|_h$ in $H_h$.
Let $I$ be the unit operator.

\par Unless otherwise stated, below $H_h$ is the space of functions given on $\bar{\omega}_h$ and equal 0 on $\partial\omega_h$, with the $L^2$-type inner product 
\[
 (v,w)_h=h_1\ldots h_n\sum\limits_{x_{\*i}\in\omega_h}v(x_{\*i})w(x_{\*i}),\ \ x_{\*i}=(i_1h_1,\ldots,i_nh_n),\ \ \*i=(i_1,\ldots,i_n).
\]
We define the simplest difference approximation $L_h:=a_1^2\Lambda_1+\ldots+a_n^2\Lambda_n$ of $L$. 
It is well-known  that
\begin{gather}
 \tfrac{4}{X_k^2}I<-\Lambda_k=-\Lambda_k^*<\tfrac{4}{h_k^2}I,\ \
 \tfrac23 I<s_{kN}=s_{kN}^*<I,
\label{Lak sNk}\\
 0<-L_h=-L_h^*<4\big(\tfrac{a_1^2}{h_1^2}+\ldots+\tfrac{a_n^2}{h_n^2}\big)I.
\label{L_h}
\end{gather}
Hereafter all the operator inequalities concern self-adjoint operators in $H_h$.

\par For completeness, we briefly recall the derivation of the method under consideration.
First we replace
%rewrite formally 
the acoustic wave equation \eqref{hyperb2eq} 
%in the form of 
with a system of equations with second order partial derivative either in $t$, or in $x_k$:
\begin{gather}
\rho(x)\partial_t^2u(x,t)-(a_1^2u_{11}(x,t)+\ldots+a_n^2u_{nn}(x,t)) =f(x,t)\ \ \text{in}\ \ Q_T,
\nonumber
\\
 u_{kk}(x,t):=\partial_k^2u(x,t),\ \ 1\leq k\leq n,\ \ \text{in}\ \ Q_T,
\label{hyperb2eq ukk}
\end{gather}
where $u_{11},\ldots,u_{nn}$ are the additional functions to seek.
%auxiliary sought functions.
Differentiating Eq. \eqref{hyperb2eq}, we get
\begin{gather}
\rho\partial_t^4u=\partial_t^2(Lu+f)=L\big[\tfrac{1}{\rho}(Lu+f)\big]+\partial_t^2f.
\label{pt2 hyperb2eq}
\end{gather}
Thus, we can sequentially write
\begin{gather} \rho\Lambda_tu=\rho\partial_t^2u+\onetwelve h_t^2\rho\partial_t^4u+\mathcal{O}(h_t^4)
\nonumber\\
 =\big(\rho I+\onetwelve h_t^2L\big)\big(\tfrac{1}{\rho}Lu\big)+f+\onetwelve h_t^2\big(\partial_t^2f+L\tfrac{f}{\rho}\big)
 +\mathcal{O}(h_t^4)
\nonumber\\
 =\big(\rho I+\onetwelve h_t^2L_h\big)\tfrac{1}{\rho}\big(a_1^2u_{11}+\ldots+a_n^2u_{nn}\big)+f_{\*h}+\mathcal{O}(|\*h|^4),\ \ \text{with}\ \ 
 f_{\*h}:=f+\onetwelve h_t^2\big(\Lambda_tf+L_h\tfrac{f}{\rho}\big).
\label{main appr err gen}
\end{gather}
Next, since the auxiliary equation \eqref{hyperb2eq ukk} is an ordinary differential equation in $x_k$, its well-known Numerov approximation leads to the formula
\begin{gather}
 s_{kN}u_{kk}-\Lambda_ku=\mathcal{O}(h_k^4)\ \ \text{on}\ \ \omega_{\*h},\ \ 1\leq k\leq n.
\label{appr error 2}
\end{gather}

\par One can omit the residues in the derived expansions \eqref{main appr err gen}-\eqref{appr error 2} and thus obtain \textit{the three-level semi-explicit in time vector compact scheme} 
\begin{gather} 
\rho\Lambda_tv-\big(\rho I+\tfrac{1}{12}h_t^2L_h\big)\tfrac{1}{\rho}(a_1^2v_{11}+\ldots+a_n^2v_{nn})=f_{\*h}
 \ \ \text{on}\ \ \omega_{\*h},
\label{fds 1}\\
 s_{kN}v_{kk}=\Lambda_kv\ \ \text{on}\ \ \omega_{\*h},\ \ 1\leq k\leq n,
\label{fds 2}
\end{gather}
with the main sought function $v\approx u$
given
%defined 
on $\bar{\omega}_h\times\overline\omega_{h_t}$
and additional ones
%and auxiliary sought functions $v\approx u$ and 
$v_{11}\approx u_{11},\ldots,v_{nn}\approx u_{nn}$ 
given
%defined 
%on $\bar{\omega}_h\times\overline\omega_{h_t}$ and 
$\bar{\omega}_h\times\omega_{h_t}$.
%, respectively.
The collection of these functions constitutes the sought vector-function.

\par The mesh
%discrete 
boundary conditions
\begin{gather} v|_{\partial\omega_{\*h}}=g,\ \
a_k^2v_{kk}|_{\partial\omega_{\*h}}=g_k,\ \ 1\leq k\leq n,
\label{fds bc}
\end{gather}
are added according to the boundary condition $u|_{\Gamma_T}=g$ and the acoustic wave equation in $\bar{Q}_T$: 
\begin{gather*}
 g_k=\rho\partial_t^2g-\sum_{1\leq l\leq n,\,l\neq k}a_l^2\partial_l^2g-f\ \ \text{for}\ \ x_k=0,X_k,\ \
\\ 
g_k=a_k^2\partial_k^2g\ \ \text{for}\ \ x_l=0,X_l,\,\ 1\leq l\leq n,\ l\neq k,
\end{gather*}
on parts of $\Gamma_T$.
In the simplest case $g=0$, 
the right-hand side of the former formula equals $-f$ and of the latter one is zero.
%the right-hand sides of these formulas equal $-f$ and $0$, respectively.

\par In addition, the function $v^m$ should be found at the first time level $m=1$ with the 4th order of accuracy.
This is often done explicitly by using Taylor's formula and the wave equation, for example, see \cite{BTT18,DZ13b,JG20,JG23,STT19}.
These cumbersome formulas use higher-order derivatives of
%(or differences) of 
%the initial functions
$u_0$ and $u_1$; they cannot be used  
%that is not satisfactory 
when $u_0$ and $u_1$
%they 
are nonsmooth.
Instead, like in \cite{Z94,ZC23,ZL24}, we apply an equation for $v^1$ similar to the above main equations \eqref{fds 1}-\eqref{fds 2} of the scheme:
\begin{gather}
 \rho(\delta_tv)^0
 =\tfrac12 h_t\big (\rho I+\tfrac{1}{12}h_t^2L_h\big)\tfrac{1}{\rho}\big(a_1^2v_{11}^0+\ldots+a_n^2v_{nn}^0\big)
 +u_{1\*h}+\tfrac{1}{2}h_tf_{\*h}^0
  \ \ \text{on}\ \ \omega_h,
\label{fds ic2 gen}\\
s_{kN}v_{kk}^0=\Lambda_kv^0\ \ \text{on}\ \ \omega_h,\ \ 1\leq k\leq n,
\label{fds ic2 vkk}
\end{gather}
with the specific $u_{1\*h}\approx \rho u_1$ and $f_{\*h}^0\approx f_0=f|_{t=0}$ such that
\begin{gather}
u_{1\*h}:=\big(\rho I+\tfrac16 h_t^2L_h\big)u_1,\ \ f_{\*h}^0:=f_{dh_t}^{(0)}+\onetwelve h_t^2L_h\tfrac{f_0}{\rho}\ \ \text{on}\ \ \omega_h,
\label{fds ic2}\\
{f}_{dh_t}^{(0)}=\tfrac{7}{12}f^0+\half f^1-\onetwelve f^2\ \ \text{or}\ \
{f}_{dh_t}^{(0)}=\tfrac13f^0+\tfrac23f^{1/2},\ \ \text{with}\ \  f^{1/2}:=f|_{t=h_t/2}.
\label{ftd02}
\end{gather}
The values of $v_{kk}^0$ on $\partial\omega_h$ can be taken as in \eqref{fds bc} (or, for smooth $u_0$, using the values of $\partial_k^2u_0$ on $\partial\Omega$).
Due to such formulas, the approximation error of Eq. \eqref{fds ic2 gen} has the 4th order
%These formulas ensure the following estimate for the approximation error of Eq. \eqref{fds ic2 gen}
\begin{gather} \hat{\psi}^0:=\rho(\delta_tu)^0-\tfrac12 h_t\big(\rho I+\tfrac{1}{12}h_t^2L_h\big)
 \tfrac{1}{\rho}\big(a_1^2u_{11\,0}+\ldots+a_n^2u_{nn\,0}\big)-u_{1\*h}-\tfrac{1}{2}h_tf_{\*h}^0
 =\mathcal{O}(|\*h|^4)
\ \ \text{on}\ \ \omega_h,
\label{hat psi 0 gen}
\end{gather}
see  \cite{ZL24}.
Hereafter $y_0:=y|_{t=0}$ for any $t$-dependent function $y$.

\par In numerical experiments below, the second (two-level) formula in \eqref{ftd02} together with more suitable forms of Eqs. \eqref{fds 1} and \eqref{fds ic2 gen} are applied 
%Notice that, in our computations below, we exploit the second (two-level) formula in \eqref{ftd02} and more practically convenient forms of Eqs. \eqref{fds 1} and \eqref{fds ic2 gen}:
\begin{gather*}
 \rho\Lambda_tv=\big(\rho I+\tfrac{1}{12}h_t^2L_h\big)\tfrac{1}{\rho}(a_1^2v_{11}+\ldots+a_n^2v_{nn}+f)+\tfrac{1}{12}h_t^2\Lambda_tf\ \ \text{on}\ \ \omega_{\*h},
\\
\rho(\delta_tv)^0=\tfrac12 h_t\big\{\rho(a_1^2v_{11}^0+\ldots+a_n^2v_{nn}^0)+{f}_{dh_t}^{(0)}
 +\tfrac{1}{12}h_t^2L_h\big[\tfrac{1}{\rho}\big(a_1^2v_{11}^0+\ldots+a_n^2v_{nn}^0+f^0\big)\big]\big\}
 +u_{1\*h}\ \ \text{on}\ \ \omega_h.
\end{gather*}

\par The described scheme is semi-explicit in time.
The function $v^{m+1}$ is found explicitly on $\omega_h$ from \eqref{fds 1} for $1\leq m\leq M-1$ or 
\eqref{fds ic2 gen} for $m=0$
when the auxiliary functions
$v_{11}^m,\ldots,v_{nn}^m$ are known on $\bar{\omega}_h$.
These functions 
are computed from the simple 1D three-point difference equations \eqref{fds 2} for $1\leq m\leq M-1$ or \eqref{fds ic2 vkk} for $m=0$ together with the boundary conditions from \eqref{fds bc} when $v^m$ is known.
Thus there exists an analogy in implementation with the ADI methods, for example, see  \cite{M90}.
On the other hand, the functions $v_{11}^m,\ldots,v_{nn}^m$ are independent and allow us multithreaded computation
%can be computed in parallel 
(that we exploit in Section \ref{numer results}), 
and we only need to store their weighted sum in order to save computer memory.
%only their weighted sum is used to compute $v^{m+1}$ that allows one to save computer memory.

\section{Stability and the 4th order error bound for the semi-explicit vector compact scheme}
\label{stab and error bounds}
\setcounter{equation}{0}
\setcounter{lemma}{0}
\setcounter{theorem}{0}
\setcounter{remark}{0}

\par We first prove an auxiliary general result.
%For functions $y$: $\{t_m\}_{m=0}^{M-1}\to H_h$, we define the norm
%\[
%\|y\|_{\tilde{L}_{h_t}^1(H_h)}:=\tfrac12 h_t\|y^0\|_h+h_t\sum_{m=1}^{M-1}\|y^m\|_h.
%\]
\begin{lemma}
Consider an abstract three-level scheme
\begin{gather}
 B_h\Lambda_tv+A_hv=\varphi\ \ \text{in}\ \ H_h\ \ \text{on}\ \ \omega_{h_t},
\label{fds 1 abs}\\
 B_h(\delta_tv)^0+\tfrac12 h_t A_hv^0=u^{(1)}+\tfrac12 h_t\varphi^0\ \ \text{in}\ \ H_h
\label{fds ic2 abs}
\end{gather}
for a function $v$: $\bar{\omega}_{h_t}\to H_h$.
Here $B_h=B_h^*>0$ and $A_h=A_h^*>0$ are 
%any 
operators acting in a Euclidean space $H_h$ 
such that
%and related by the inequality
\begin{gather}
 \tfrac14 h_t^2A_h\leq (1-\ve_0^2)B_h\ \ \text{for some}\ \ 0<\ve_0<1.
\label{stab cond 2 gen}
\end{gather}
Then the following stability bounds  hold
\begin{gather}
 \max_{0\leq m\leq M}
 \max\big\{\ve_0\|v^m\|_{B_h},\,\|I_{h_t}^m\bar{s}_tv\|_{A_h}\big\}
\leq\|v^0\|_{B_h}+2\|A_h^{-1/2}u^{(1)}\|_h
+2\big\|A_h^{-1/2}\vp\big\|_{\tilde{L}_{h_t}^1(H_h)},
\label{gen stab bound 0}\\
 \ve_0\max_{1\leq m\leq M}\|A_h^{-1/2}B_h\bar{\delta}_tv\|_h
\leq (1+\ve_0)\|v^0\|_{B_h}+(3+2\ve_0)\big(\|A_h^{-1/2}u^{(1)}\|_h
+\big\|A_h^{-1/2}\vp\big\|_{\tilde{L}_{h_t}^1(H_h)}\big). 
 \label{gen stab bound 0 plus}
\end{gather}
Here, we use the norm
$
\|y\|_{\tilde{L}_{h_t}^1(H_h)}:=\tfrac12 h_t\|y^0\|_h+h_t\sum_{m=1}^{M-1}\|y^m\|_h\ \ \text{for functions}\ \ y: \{t_m\}_{m=0}^{M-1}\to H_h.
$
\end{lemma}
\begin{proof}
Bound \eqref{gen stab bound 0} follows from the known bound in the weak energy norm
\begin{gather*}
 \max_{0\leq m\leq M}
 \max\big\{\big[\|v^m\|_{B_h}^2-\tfrac14h_t^2\|v^m\|_{A_h}^2\big]^{1/2},\,\|I_{h_t}^m\bar{s}_tv\|_{A_h}\big\}
\\
\leq\big[\|v^0\|_{B_h}^2-\tfrac14h_t^2\|v^0\|_{A_h}^2\big]^{1/2}
+2\|A_h^{-1/2}u^{(1)}\|_h
+2\big\|A_h^{-1/2}\vp\big\|_{\tilde{L}_{h_t}^1(H_h)},
\end{gather*}
see \cite[Theorem 1]{ZK20} and \cite[Theorem 2]{ZCh20} (for the weight $\sigma=0$), and from condition \eqref{stab cond 2 gen}.
\par To derive bound \eqref{gen stab bound 0 plus},
we apply the operator $I_{h_t}$ to Eq. \eqref{fds 1 abs} and get 
\[
B_h(\delta_tv^m-\delta_tv^0)+A_hI_{h_t}^mv=I_{h_t}^m\varphi,\ \ 0\leq m\leq M-1.
\]
Due to the initial condition \eqref{fds ic2 abs} and the elementary formula
\begin{equation}
\tfrac12h_tv^0+I_{h_t}^mv
=I_{h_t}^m\bar{s}_tv+\tfrac12h_tv^m,\ \ 0\leq m\leq M,
\label{elem form}
\end{equation}
we further obtain
\[
B_h\delta_tv^m=-A_h(I_{h_t}^m\bar{s}_tv+\tfrac12h_tv^m)+u^{(1)}+\tfrac12 h_t\varphi^0+I_{h_t}^m\varphi,\ \ 0\leq m\leq M-1.
\]
Next, we apply the operator $A_h^{-1/2}$ to both sides of this equality, take the norm in $H_h$ and, using condition \eqref{stab cond 2 gen}, derive
\begin{gather*}
\|A_h^{-1/2}B_h\bar{\delta}_tv^{m+1}\|_h
\leq \|I_{h_t}^m\bar{s}_tv\|_{A_h}+\tfrac12h_t\|v^m\|_{A_h}+\|A_h^{-1/2}u^{(1)}\|_h+\tfrac12 h_t\|A_h^{-1/2}\varphi^0\|_h
+I_{h_t}^m\|A_h^{-1/2}\varphi\|_h
\\
\leq (1+\ve_0^{-1})\max\big\{\ve_0\|v^m\|_{B_h},\,\|I_{h_t}^m\bar{s}_tv\|_{A_h}\big\}
+\|A_h^{-1/2}u^{(1)}\|_h
+\big\|A_h^{-1/2}\vp\big\|_{\tilde{L}_{h_t}^1(H_h)}.
\end{gather*}
Now the first bound of this Lemma implies the second one. 
\end{proof}
\begin{remark}
\label{rem on Iv}
Due to formula \eqref{elem form} and condition \eqref{stab cond 2 gen}, we obtain
\[
\|I_{h_t}^mv\|_{A_h}
\leq\tfrac12h_t\|v^0\|_{A_h}+\|I_{h_t}^m\bar{s}_tv\|_{A_h}+\tfrac12h_t\|v^m\|_{A_h}
\leq\|v^0\|_{B_h}+\|I_{h_t}^m\bar{s}_tv\|_{A_h}+\|v^m\|_{B_h},\ \ 0\leq m\leq M.
\]
Thus, an additional bound holds 
\[
\ve_0\max_{0\leq m\leq M}
\|I_{h_t}^mv\|_{A_h}
\\
\leq (2+2\ve_0)\big(\|v^0\|_{B_h}
+\|A_h^{-1/2}u^{(1)}\|_h
+\big\|A_h^{-1/2}\vp\big\|_{\tilde{L}_{h_t}^1(H_h)}\big).
\]
\end{remark}
\par We consider the inhomogeneous version of Eqs. \eqref{fds 2} and \eqref{fds ic2 vkk}:
\begin{gather}
s_{kN}v_{kk}^m=\Lambda_kv^m+b_k^m\ \ \text{on}\ \ \omega_h,\ \ 0\leq m\leq M-1,\ \ 1\leq k\leq n,
\label{fds 2m}
\end{gather}
with the given functions $b_1,\ldots,b_n$.
In practical computations, one never has $b_1=\ldots=b_n=0$ due to the presence of the round-off errors.
Therefore, an impact of $b_1,\ldots,b_n$ on $v$ has to be studied.
Importantly, one must do that to prove error bounds for the scheme.
%, so the influence of  $b_1,\ldots,b_n$ on the solution needs to be investigated.
%Importantly, this is also required to prove any error bound.

\par We define the following two self-adjoint operators in $H_h$:
%To state the stability theorem, we introduce the following self-adjoint operators in $H_h$:
\[
 E_h:=-(a_1^2s_{1N}^{-1}\Lambda_1+\ldots+a_n^2s_{nN}^{-1}\Lambda_n),\ \  I_{\rho\*h}:=\tfrac{1}{\rho}I+\tfrac{1}{12}h_t^2\big(\tfrac{1}{\rho}L_h\big)\tfrac{1}{\rho}I.
\]
The second inequality \eqref{Lak sNk} can be rewritten as  $I<s_{1N}^{-1}<(3/2)I$ and thus 
$-L_h<E_h<-(3/2) L_h$.
Clearly also
$I_{\rho\*h}=I_{\rho\*h}^*
<\rho^{-1}I\leq\urho^{-1}I$.

\par We introduce the 
first CFL-type condition on $h_t$ and $h_1,\ldots,h_n$:
%following first CFL-type condition on the steps in time and space
\begin{gather}
 \tfrac{1}{3}h_t^2\big(\tfrac{a_1^2}{h_1^2}+\ldots+\tfrac{a_n^2}{h_n^2}\big)\leq(1-\ve)\urho
 \ \ \text{for some}\ \ 0\leq\ve<1.
\label{stab cond 1 gen}
\end{gather}
Then with the help of inequality \eqref{L_h} we have
\[
-\big(\tfrac{1}{12}h_t^2\big(\tfrac{1}{\rho}L_h\big)\tfrac{w}{\rho},w\big)_h
=-\tfrac{1}{12}h_t^2\big(L_h\tfrac{w}{\rho},\tfrac{w}{\rho}\big)_h
<(1-\ve)\urho\big(\tfrac{w}{\rho},\tfrac{w}{\rho}\big)_h
\leq (1-\ve)\big(\tfrac{1}{\rho}w,w\big)_h\ \ \text{for any}\ \ w\in H_h
\]
and, consequently, we can also bound $I_{\rho\*h}$ from below and get the two-sided bounds
\begin{gather}
\ve \tfrac{1}{\bar{\rho}}I
\leq\ve \tfrac{1}{\rho}I
<I_{\rho\*h}<\tfrac{1}{\rho}I\leq\tfrac{1}{\urho}I.
\label{lower bound for I rho h}
\end{gather}
For $0<\ve<1$, these bounds are equivalent to the following ones
\begin{equation}
\urho I
\leq \rho I
<I_{\rho\*h}^{-1}<\ve^{-1}\rho I\leq\ve^{-1}\bar{\rho}I.
\label{bounds for inv I rho h} 
\end{equation}
\par The new stability result on the semi-explicit vector compact scheme is as follows.
\begin{theorem}
\label{theo:stab bound 2}
%Let 
Assume that the stability conditions \eqref{stab cond 1 gen} and
\begin{gather}
 \tfrac14 h_t^2E_h
 \leq (1-\ve_0^2)\urho I\ \ \text{for some}\ \ 0<\ve_0<1
\label{stab cond 2 gener}
\end{gather}
are valid, and $g=g_1=\ldots=g_n=0$ in the boundary conditions \eqref{fds bc}. Then 
%be valid.
%Then, for scheme \eqref{fds 1}, \eqref{fds 2m} and \eqref{fds ic2 gen} under the boundary conditions \eqref{fds bc} with $g=g_1=\ldots=g_n=0$,
the following stability bounds hold
\begin{gather}
 \max_{0\leq m\leq M}
 \max\big\{\ve_0\|v^m\|_{E_h},\,\|I_{h_t}^m\bar{s}_tE_hv\|_{I_{\rho\*h}}\big\}
\leq\|v^0\|_{E_h}+2\big\|\tfrac{1}{\rho}u_{1\*h}\big\|_{I_{\rho\*h}^{-1}}
+2\big\|I_{\rho\*h}^{-1/2}\big[\tfrac{1}{\rho}(f_{\*h}+\beta_{\*h})\big]\big\|_{\tilde{L}_{h_t}^1(H_h)},
\label{stab bound 0}\\
 \ve_0\max_{1\leq m\leq M}\|\bar{\delta}_tv^m\|_{I_{\rho\*h}^{-1}}
\leq (1+\ve_0)\|v^0\|_{E_h}+(3+2\ve_0)\big(\big\|\tfrac{1}{\rho}u_{1\*h}\big\|_{I_{\rho\*h}^{-1}}
+\big\|I_{\rho\*h}^{-1/2}\big[\tfrac{1}{\rho}(f_{\*h}+\beta_{\*h})\big]\big\|_{\tilde{L}_{h_t}^1(H_h)}\big), 
 \label{stab bound 0 plus}
\end{gather}
for any functions $f_{\*h},b_1,\ldots,b_n$: $\{t_m\}_{m=0}^{M-1}\to H_h$ and $v^0,u_{1\*h}\in H_h$  (thus, $f_{\*h}$ and $u_{1\*h}$ are not only those specific functions defined \eqref{fds ic2}-\eqref{ftd02}) and
\begin{gather}
\beta_{\*h}^m:=\big(\rho I+\tfrac{1}{12}h_t^2L_h\big)\tfrac{1}{\rho}(a_1^2s_{1N}^{-1}b_1^m+\ldots+a_n^2s_{nN}^{-1}b_n^m)\ \ \text{in}\ \ H_h,\ \ 0\leq m\leq M-1.
\label{fds 1 sep beta gen}
\end{gather}
\end{theorem}
\begin{proof}
Since $g=g_1=\ldots=g_n=0$,
clearly
$v$: $\overline\omega_{h_t}\to H_h$ and $v_{kk}$: $\{t_m\}_{m=0}^{M-1}\to H_h$, $1\leq k\leq n$.
Acting as in \cite{ZL24}, we can express $v_{kk}$ from Eq. \eqref{fds 2m}:
\begin{gather}
 v_{kk}^m=s_{kN}^{-1}(\Lambda_kv^m+b_k^m)\ \ \text{in}\ \ H_h,\ \ 
 1\leq k\leq n,\ \ 0\leq m\leq M-1.
\label{vkk thr v}
\end{gather}
Using these expressions in Eqs. \eqref{fds 1} and \eqref{fds ic2 gen} leads to the equations for~$v$ only:
\begin{gather*}
\rho\Lambda_tv^m+\big(\rho I+\tfrac{1}{12}h_t^2L_h\big)\tfrac{1}{\rho}E_hv^m
  =f_{\*h}^m+\beta_{\*h}^m
\ \ \text{in}\ \ H_h,
  \ \ 1\leq m\leq M-1,
\label{fds 1 sep gen}\\
 \rho(\delta_tv)^0+\tfrac12 h_t\big(\rho I+\tfrac{1}{12}h_t^2L_h\big)\tfrac{1}{\rho}E_hv^0
 =u_{1\*h}+\tfrac{1}{2}h_t(f_{\*h}^0+\beta_{\*h}^0)
\ \ \text{in}\ \ H_h;
\label{fds 2 sep gen}
\end{gather*}
recall that $\beta_{\*h}$ has been introduced in \eqref{fds 1 sep beta gen}.

\par Importantly, both the equations can be symmetrized by applying the operator $E_h\big(\rho^{-1}I\big)$ to them
\begin{gather} E_h\Lambda_tv^m+E_hI_{\rho\*h}E_hv^m
  =E_h\big(\tfrac{1}{\rho}(f_{\*h}+\beta_{\*h})^m\big)
\ \ \text{in}\ \ H_h,
  \ \ 1\leq m\leq M-1,
\label{fds 1 sep sym}\\
 E_h(\delta_tv)^0+\tfrac12 h_t E_hI_{\rho\*h}E_hv^0
 =E_h\big(\tfrac{1}{\rho}u_{1\*h}\big)+\tfrac{1}{2}h_t E_h\big(\tfrac{1}{\rho}(f_{\*h}^0+\beta_{\*h}^0)\big)
\ \ \text{in}\ \ H_h.
\label{fds 2 sep sym}
\end{gather}
Here $E_hI_{\rho\*h}E_h=(E_hI_{\rho\*h}E_h)^*>0$ under condition \eqref{stab cond 1 gen} due to inequalities \eqref{lower bound for I rho h}, therefore, these equations can be treated as a particular case of scheme \eqref{fds 1 abs}-\eqref{fds ic2 abs} for $B_h=E_h$ and $A_h=E_hI_{\rho\*h}E_h$.

\par Now the general stability condition \eqref{stab cond 2 gen} looks as
\[
 \tfrac{1}{4}h_t^2 E_hI_{\rho\*h}E_h\leq (1-\ve_0^2)E_h\ \ \text{for some}\ \ 0<\ve_0<1.
\]
Recall that $I_{\rho\*h}<\urho^{-1}I$, thus the last condition follows from $(1/4)h_t^2 E_h^2\leq (1-\ve_0^2)\urho E_h$ or, equivalently, from \eqref{stab cond 2 gener}.

\par Next, for any $w\in H_h$, we have
\begin{gather*}
\|w\|_{A_h}^2=(E_hI_{\rho\*h}E_hw,w)_h
=\|E_hw\|_{I_{\rho\*h}}^2,
\\
\|A_h^{-1/2}B_hw\|_h^2=(A_h^{-1}B_hw,B_hw)_h
=\big((E_hI_{\rho\*h}E_h)^{-1}E_hw,E_hw\big)_h
=(I_{\rho\*h}^{-1}w,w)_h=\|w\|_{I_{\rho\*h}^{-1}}^2.
\end{gather*}
Consequently, the general stability bounds \eqref{gen stab bound 0}-\eqref{gen stab bound 0 plus} take the form of the stated bounds \eqref{stab bound 0}-\eqref{stab bound 0 plus}.
\end{proof}

\par Since $E_h<\frac32(-L_h)$ in $H_h$,
both
the stability conditions \eqref{stab cond 1 gen} and \eqref{stab cond 2 gener} follow from the inequality
\begin{gather*}
h_t^2\big(\tfrac{a_1^2}{h_1^2}+\ldots+\tfrac{a_n^2}{h_n^2}\big)\leq \min\{3(1-\ve),\tfrac23(1-\ve_0^2)\}\urho.
\label{stab cond 1 gen expl}
\end{gather*}

\par The right-hand sides of bounds \eqref{stab bound 0}-\eqref{stab bound 0 plus} can be simplified. Due to the formula
$\rho^{-1}\beta_{\*h}
=I_{\rho\*h}(a_1^2s_{1N}^{-1}b_1+\ldots+a_n^2s_{nN}^{-1}b_n)$ and the operator inequalities 
\eqref{lower bound for I rho h}, \eqref{bounds for inv I rho h} and $s_{kN}^{-1}\leq (3/2) I$, for $0<\ve<1$, we get
\begin{gather}
\big\|I_{\rho\*h}^{-1/2}\big[\tfrac{1}{\rho}(f_{\*h}+\beta_{\*h})\big]\big\|_{\tilde{L}_{h_t}^1(H_h)}
\leq \ve^{-1/2}\big\|\tfrac{f_{\*h}}{\sqrt{\rho}}\big\|_{\tilde{L}_{h_t}^1(H_h)}
+\big\|I_{\rho\*h}^{1/2}(a_1^2s_{1N}^{-1}b_1+\ldots+a_n^2s_{nN}^{-1}b_n)\big\|_{\tilde{L}_{h_t}^1(H_h)}
\\
\leq\urho^{-1/2}\Big(\ve^{-1/2}\|f_{\*h}\|_{\tilde{L}_{h_t}^1(H_h)}
+\tfrac32a_{\max}^2\sum_{k=1}^n\|b_k\|_{\tilde{L}_{h_t}^1(H_h)}\Big),
\label{est f beta}\\
\big\|\tfrac{1}{\rho}u_{1\*h}\big\|_{I_{\rho\*h}^{-1}}=\big\|I_{\rho\*h}^{-1/2}\big(\tfrac{1}{\rho}u_{1\*h}\big)\big\|_h\leq (\ve\urho)^{-1/2}\|u_{1\*h}\|_h,\ \
\urho^{1/2}\|\bar{\delta}_tv^m\|_h\leq\|\bar{\delta}_tv^m\|_{I_{\rho\*h}^{-1}}.
\label{aux bounds}
\end{gather}

\par Let $1\leq k\leq n$. For functions $w$ defined on $\bar{\omega}_{hk}$, we introduce a difference operator and seminorms
\begin{gather*} (\bar{\delta}_kw)_i=\tfrac{1}{h_k}(w_i-w_{i-1}),\,\
1\leq i\leq N_k,
\ \
\|w\|_{k}=\Big(h_k\sum_{i=1}^{N_k-1}w_i^2\Big)^{1/2},\ \ \|w\|_{k*}=\Big(h_k\sum_{i=1}^{N_k}w_i^2\Big)^{1/2}.
\end{gather*}

\par For functions $w$ 
given on
%defined on 
$\bar{\omega}_h$, we have
%recall the seminorm
$\|w\|_h=\|\ldots\|w\|_1\ldots\|_n$.
%and 
Introduce
%define 
seminorms
$\|w\|_{h,k*}$  which defer from $\|w\|_h$ by replacing $\|w\|_k$ with $\|w\|_{k*}$ and then
a mesh counterpart of the squared norm in the Sobolev subspace $H_0^1(\Omega)$:
\[
 \|w\|_{H_h^1}^2:=\|(-L_h)^{1/2}w\|_h^2=(-L_hw,w)_h=\sum_{k=1}^na_k^2\|\bar{\delta}_kw\|_{h,k*}^2.
\]
Note the two-sided bounds
$\|w\|_{H_h^1}^2\leq \|w\|_{E_h}^2\leq \tfrac32\|w\|_{H_h^1}^2.$
We also define the mesh energy norm 
$\|y\|_{\mathcal{E}_h}:=\big(\|\bar{\delta}_ty\|_h^2+\|y\|_{H_h^1}^2\big)^{1/2}$
that is needed
%used below
in Section \ref{numer results}.

\par The next theorem states
%The next main result concerns
the 4th order error bound in the enlarged mesh energy norm, and it is derived from Theorem~\ref{theo:stab bound 2}.
%Notice that 
In this theorem, $g$ and $f|_{\Gamma_T}$ can be arbitrary (not only zero).
%in it.
\begin{theorem}
\label{theo:error bound}
Let
%Under 
the stability conditions \eqref{stab cond 1 gen} with $0<\ve<1$ and \eqref{stab cond 2 gener}
hold, and $v^0=u^0$ on $\bar{\omega}_h$.
Then the 
%following 
4th order error bound in the enlarged mesh energy norm for scheme
\eqref{fds 1}-\eqref{fds ic2 vkk} 
%with $v^0=u^0$ on $\bar{\omega}_h$ 
holds
\begin{gather}
 \sqrt{\ve}\ve_0\max_{1\leq m\leq M}\big(\|\bar{\delta}_t(u-v)^m\|_h
 +\|(u-v)^m\|_{H_h^1}
 +\sqrt{\ve}\|I_{h_t}^mL_h(u-v)\|_h\big)=\mathcal{O}(|\*h|^4).
\label{error bound}
\end{gather}
Here and below the $\mathcal{O}(|\*h|^4)$-terms do not depend on $\ve$ and $\ve_0$.
%Hereafter the $\mathcal{O}(|\*h|^4)$-terms are independent of $\ve$ and $\ve_0$.
\end{theorem}
\begin{proof}
First, the
%The
approximation errors of Eqs. \eqref{fds 1}, \eqref{fds 2} and \eqref{fds ic2 vkk} are expressed by the formulas
\begin{gather*}
 \psi:=\rho\Lambda_tu
 -\big(\rho I+\onetwelve h_t^2L_h\big)\tfrac{1}{\rho}\big(a_1^2u_{11}+\ldots+a_n^2u_{nn}\big)+f_{\*h}
 \ \ \text{on}\ \ \omega_{\*h},
\\
\psi_{kk}:=s_{kN}u_{kk}-\Lambda_ku
\ \ \text{on}\ \ \omega_{\*h},\ \ 
\psi_{kk}^0:=s_{kN}\partial_k^2u_0-\Lambda_ku_0
\ \ \text{on}\ \ \omega_h,\ \ 1\leq k\leq n,
\end{gather*}
whereas the approximation error of Eq. \eqref{fds ic2 gen} has been given by \eqref{hat psi 0 gen}.
According to formulas \eqref{main appr err gen}, \eqref{appr error 2} and \eqref{hat psi 0 gen}, all the approximation errors are of the 4th order in the mesh uniform norm
\begin{gather} \max_{\omega_{\*h}}|\psi|+\max_{\omega_h}|\hat{\psi}^0|
 +\max_{0\leq m\leq M-1}\max_{\omega_h}\big(|\psi_{11}^m|+\ldots+|\psi_{nn}^m|\big)
 =\mathcal{O}(|\*h|^4).
\label{appr error bound 2}
\end{gather}

\par Due to the equations for $v$ and $v_{11},\ldots,v_{nn}$ as well as the definitions of $\psi$, $\psi_{kk}$ and $\hat{\psi}^0$,
the errors $r=u-v$ and $r_{11}=u_{11}-v_{11},\ldots,r_{nn}=u_{nn}-v_{nn}$ satisfy the following system of equations
\begin{gather}
 \rho\Lambda_tr
 -\big(\rho I+\onetwelve h_t^2L_h\big)\tfrac{1}{\rho}\big(a_1^2r_{11}+\ldots+a_n^2r_{nn}\big)=\psi
 \ \ \text{on}\ \ \omega_{\*h},
\nonumber\\
s_{kN}r_{kk}-\Lambda_kr=\psi_{kk}\ \ \text{on}\ \ \omega_{\*h},\ \
s_{kN}r_{kk}^0-\Lambda_kr^0=\psi_{kk}^0\ \ \text{on}\ \ \omega_h,\ \ 1\leq k\leq n,
\label{eq for rkk}\\
 \rho(\delta_tr)^0
 -\tfrac12 h_t\big (\rho I+\tfrac{1}{12}h_t^2L_h\big)\tfrac{1}{\rho}\big(a_1^2r_{11}^0+\ldots+a_n^2r_{nn}^0\big)=\hat{\psi}^0\ \ \text{on}\ \ \omega_h,
 \nonumber
\end{gather}
with the approximation errors on the right.
The respective boundary and initial data are zero:
\[ r|_{\partial\omega_{\*h}}=0,\ \ r_{kk}|_{\partial\omega_{\*h}}=0,\ \ 1\leq k\leq n,\ \ r^0=0.
\]
Note that here $\psi$ (with $\psi^0:=0$), $\psi_{kk}$ and $\hat{\psi}^0$ play the role of $f_{\*h}$, $b_k$ and $u_{1\*h}$, respectively, in the original scheme equations.

\par The stability bounds of Theorem~\ref{theo:stab bound 2} applied to these equations, together with Remark \ref{rem on Iv} and estimates \eqref{est f beta}, \eqref{aux bounds} and \eqref{appr error bound 2},
lead to the error bound
\begin{gather*}
\ve_0\max_{1\leq m\leq M}\max\big\{
\urho^{1/2}\|\bar{\delta}_tr^m\|_h,\,\|r^m\|_{E_h},\,
\|I_{h_t}^mE_hr\|_{I_{\rho\*h}}\big\}
\\
 \leq (3+2\ve_0) (\ve\urho)^{-1/2}\Big[\|\hat{\psi}^0\|_h
 +h_t\sum_{m=1}^{M-1}\|\psi^m\|_h
 +\tfrac32a_{\max}^2 h_t\sum_{m=0}^{M-1}\big(\|\psi_{11}^m\|_h+\ldots+\|\psi_{nn}^m\|_h\big)\Big]
 =\mathcal{O}(|\*h|^4).
\end{gather*}
Using first inequality \eqref{lower bound for I rho h} and then
$-L_h\leq E_h$ in the last term on the left, we get the error bound \eqref{error bound}.
\end{proof}

\par According to \cite{ZL24}, the derived error bound implies the additional error bound for the auxiliary unknowns in weaker ``negative'' norms.
Indeed, from Eqs. \eqref{eq for rkk} we get $r_{kk}=s_{kN}^{-1}(\Lambda_kr+\psi_{kk})$ and thus
\begin{gather*}
\|(-\Lambda_k)^{-1/2}r_{kk}\|_h
\leq\tfrac32\big(\|(-\Lambda_k)^{1/2}r\|_h+\tfrac{X_k}{2}\|\psi_{kk}\|_h\big),\ \ 1\leq k\leq n,
\end{gather*}
see \eqref{Lak sNk},
that leads to the error bound
\begin{gather*}
\sqrt{\ve}\ve_0\max_{1\leq k\leq n}\max_{0\leq m\leq M-1}\|(-\Lambda_k)^{-1/2}(\partial_k^2u^m-v_{kk}^m)\}\|_h
=\mathcal{O}(|\*h|^4).
\end{gather*}

\section{Numerical results}
\label{numer results}
\setcounter{equation}{0}
\setcounter{lemma}{0}
\setcounter{theorem}{0}
\setcounter{remark}{0}

This section is devoted to results of our 3D numerical experiments on various tests.
%They have been performed on the computer with AMD\textsuperscript{\textcopyright} Ryzen\textsuperscript{\textcopyright} 5 7600X 6-Core Processor and 64GB RAM.
The scheme code is written in Pure C language for x64-bit architecture. The figures in this section are plotted with the use of Python standard libraries.
The computer with AMD\textsuperscript{\textcopyright} Ryzen\textsuperscript{\textcopyright} 5 7600X 6-Core Processor and 64GB RAM is applied.
%The algorithm itself and the GUI have been implemented using \textit{C} and \textit{Python} languages, respectively.

\par First, for $h_1=\ldots=h_n=h$, we consider the error of some scheme in a norm and assume that its convergence rate is some $p>0$ with respect to both $h$ and $h_t$:
\[
e(h,h_t)=c_1h^p+c_2h_t^p,
\]
for given data, with $c_1>0$ and $c_2>0$. 
Recall that $p$ depends not only on the approximation order of the scheme but on the smoothness order of the data as well, and these orders are the same only for sufficiently smooth data.
Then, for some $q>1$, we get
\[
e(h/q,h_t/q)=e(h,h_t)/q^p
\]
and, thus, we can find the unknown $p$ due to the well-known Runge-type formula
\begin{equation} 
p=p(e)=\Big.\ln\frac{e(h,h_t)}{e(h/q,h_t/q)}\Big/\ln q.
\label{pr_err_order}
\end{equation}
Of course, this practical formula is only approximate. But one can expect that it becomes more and more accurate as $h$ and $h_t$ are chosen smaller and smaller.
In Examples 1 and 2, we will analyze the errors $e_{L_h^2}$ in the $H_h=L_h^2$-norm, $e_{H_h^1}$ in the $H_h^1$-norm and $e_{\mathcal{E}_h}$ in the $\mathcal{E}_h$-norm all taken at the final time moment $t = T$ (unless otherwise is stated). Below we use $q=5/3$ less than the most common $q=2$ to cut computational costs in half.
\par In all the three examples, we choose $a_1=a_2=a_3=a=\rm{const}$ even in the case of the variable speed of sound in contrast to \cite{JG23}. Also we take the cubic domains $\Omega$, i.e., $X_1=X_2=X_3=X$, and the simplest cubic mesh with $N_1=N_2=N_3=N$, i.e., $h_1=h_2=h_3=h=X/N$.
As usual, we denote $(x_1,x_2,x_3)=(x,y,z)$.

\smallskip\par{\bf Example 1.}
We take $a=1/\sqrt{3}$, $\Omega=(0,1)^3$ and $T=0.3$ 
as well as
%and the data $f=0$, $u_0$, $u_1$ and $g$ 
%defined according to 
the exact travelling wave smooth solution
\[
 u(x,y,z,t)=\cos(t-x-y-z).
\]
%that is a travelling wave. 
The data $f=0$, $u_0$, $u_1$ and $g$ are taken respectively, and all of them
%Note that all the corresponding data 
including $g$ are non-zero.

\par We select values  $N=81,135,225,375$ and $M=27,45,75,125$ so that the ratio of the consecutive values of both $N$ and $M$ is constant and equals $q=5/3$.
Note that then $\sqrt{3}a(\tau/h)=0.9<1$ ensuring stability of the scheme.

\par (a) Let first $\rho=1$ be constant. The results are collected in Table \ref{tab: tr wave solut}, where all the errors are excellent (very small, namely, less than $1.62\times 10^{-10}$ even for the least values $N=81$ and $M=27$), and the rates $p_{L_h^2}$ are very close to the approximation order 4.
The first and second values of the rates $p_{H_h^1}$ and $p_{\mathcal{E}_h}$ are also very close to 4, but their third values become less close to 4 due to the influence of round-off errors.
The included theoretical error orders 4.000 hereafter formally correspond to the values $N=M=*$; here all of them equal 4.
Notice that, in this and other Tables, always $e_{L_h^2}<e_{\mathcal{E}_h}<e_{H_h^1}$.
In this Table, the included value $k=4$ is the approximation order of the scheme.
\begin{table}[ht]
\begin{center}
\caption{Example 1(a). Errors $e_{L_h^2}$,  $e_{H_h^1}$ and $e_{\mathcal{E}_h}$ and convergence rates $p_{L_h^2}$, $p_{H_h^1}$ and $p_{\mathcal{E}_h}$, for the travelling wave solution and constant $\rho$}
\label{tab: tr wave solut}
\begin{tabular}{crrcccccc}
\hline\noalign{\smallskip}
& $N$ & $M$ & 
$e_{L_h^2}$ & $e_{H_h^1}$ & $e_{\mathcal{E}_h}$ &
$p_{L_h^2}$ & $p_{H_h^1}$ & $p_{\mathcal{E}_h}$ \\
\noalign{\smallskip}
\hline
\noalign{\smallskip}
$k=4$ & 81 & 27& 2.434899E-11 &1.618170E-10 &1.166171E-10 &$-$   &$-$    &$-$ \\
      & 135& 45& 3.186161E-12 &2.119400E-11 &1.528949E-11 &3.981 &3.979 &3.977 \\
      & 225& 75& 4.153367E-13 &2.766222E-12 &1.996804E-12 &3.989 &3.986 &3.985 \\
      & 375&125& 5.400149E-14 &4.070984E-13 &3.056195E-13 &3.994 &3.751 &3.674 \\
 & $*$&  $*$  &$-$   &$-$   &$-$                &4.000 &4.000 &4.000 \\
\hline
\end{tabular}
\end{center}
\end{table}

\par (b) Next, we take variable 
$\rho(x)=1+\sin^22\pi x\cdot\sin^22\pi y\cdot\sin^22\pi z$, with $\urho=1$.
This does not lead to a significant impact on the results that are put in
Table \ref{tab: tr wave solut var rho}. 
All the results are very close to the corresponding ones in Table~\ref{tab: tr wave solut}. 
For variable $\rho$, all the errors are even slightly less, but almost all the convergence rates are very slightly less too. 

This Table includes the additional column containing the ratio $CPU_{rel}$ of the CPU times for consecutive values of $N$ and $M$. They are close to the theoretical ratio $(5/3)^4\approx 7.72$ that is listed last in this column. Note that, for the usual $q=2$, this ratio is about twice as large: $2^4=16$.
\begin{table}[ht]
\begin{center}
\caption{Example 1(b). Errors $e_{L_h^2}$,  $e_{H_h^1}$ and $e_{\mathcal{E}_h}$, convergence rates $p_{L_h^2}$, $p_{H_h^1}$ and $p_{\mathcal{E}_h}$ and the ratio of the CPU times for the travelling wave solution and variable $\rho$}
\label{tab: tr wave solut var rho}
\begin{tabular}{crrccccccc}
\hline\noalign{\smallskip}
& $N$ & $M$ & 
$e_{L_h^2}$ & $e_{H_h^1}$ & $e_{\mathcal{E}_h}$ &
$p_{L_h^2}$ & $p_{H_h^1}$ & $p_{\mathcal{E}_h}$ & $CPU_{rel}$\\
\noalign{\smallskip}
\hline
\noalign{\smallskip}
$k=4$ & 81 & 27& 2.083224E-11 &1.405375E-10 &1.035755E-10 &$-$   &$-$   &$-$   &$-$\\
      & 135& 45& 2.725370E-12 &1.845130E-11 &1.361021E-11 &3.982 &3.974 &3.973 &7.40\\
      & 225& 75& 3.552248E-13 &2.412085E-12 &1.780023E-12 &3.989 &3.983 &3.982 &7.49\\
      & 375&125& 4.618912E-14 &3.684124E-13 &2.813984E-13 &3.993 &3.678 &3.611 &7.87\\
      & $*$&$*$&        $-$   &$-$          &$-$          &4.000 &4.000 &4.000 &7.72\\
\hline
\end{tabular}
\end{center}
\end{table}

\par
\begin{figure}[tbh!]
\begin{minipage}{0.5\textwidth}
\center{\includegraphics[width=1\linewidth]{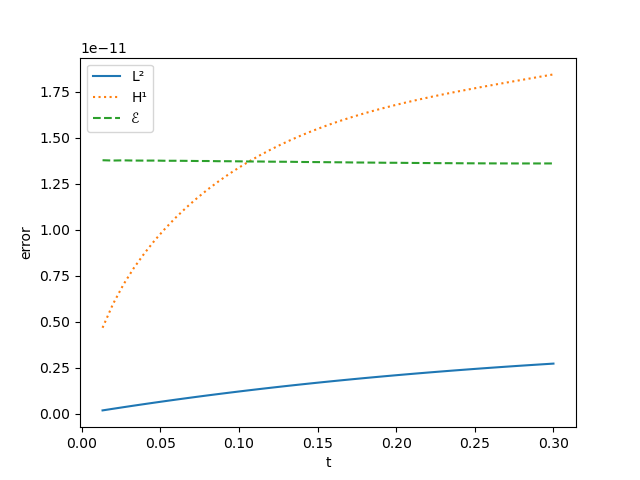}} (a) \\
 \end{minipage}
 \hfill
 \begin{minipage}{0.5\textwidth}
\center{\includegraphics[width=1\linewidth]{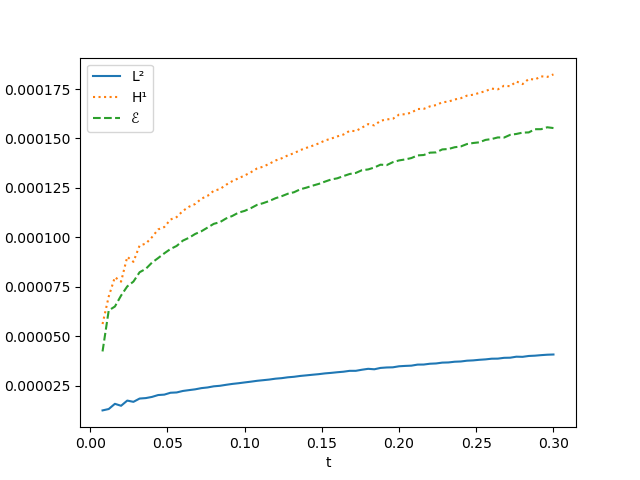}} (b) \\
\end{minipage}
\caption{Dynamics of the errors $e_{L_h^2}$,  $e_{H_h^1}$ and $e_{\mathcal{E}_h}$ in time for: (a) Example 1(b) for $N=135$ and $M=45$, (b) Example 2 in the case $u_1=w_1$ for $N=225$ and $M=75$. In case (b), $100e_{L_h^2}$ is presented instead of $e_{L_h^2}$.}
\label{fig:errors_time}
\end{figure}

 In Fig. \ref{fig:errors_time}(a), we also present dynamics of the errors $e_{L_h^2}$,  $e_{H_h^1}$ and $e_{\mathcal{E}_h}$ in time.
We observe slow and almost linear growth of $e_{L_h^2}$, much more rapid but slower than linear growth of $e_{H_h^1}$ and non-standard very slow linear decrease of  $e_{\mathcal{E}_h}$.

\smallskip\par{\bf Example 2.}
To analyze the practical error orders for the nonsmooth data, we take the following spherically symmetric initial functions
$u_0=w_1,w_2$, $ u_1=w_0,w_1$ and free term $f=w_0,w_1$
using piecewise polynomial functions with a finite support in a domain $\Omega$ symmetric with respect to the origin:
\[
w_0(r)=\begin{cases}
1,\ r<r_0\\[2mm]
0,\ r>r_0
\end{cases},\ \
w_1(r)=\begin{cases}
\dfrac{r_0-r}{r_0},\ r\leq r_0\\[2mm]
0,\ r>r_0
\end{cases},\ \
w_2(r)=\begin{cases}
\Big(\dfrac{r}{r_0}\Big)^2\Big(\dfrac{r_0-r}{r_0}\Big)^2,\ r\leq r_0\\[2mm]
0,\ r>r_0
\end{cases},
\]
with $r=|x|$ and the parameter $r_0>0$.
Clearly $w_0\in L^\infty(\mathbb{R}^+)$, $w_1\in W^{1,\infty}(\mathbb{R}^+)$ and $w_2\in W^{2,\infty}(\mathbb{R}^+)$, where  the usual Lebesgue and Sobolev spaces on $\mathbb{R}^+=(0,+\infty)$ are used.
More importantly below, since $w_0,w_1'$ and $w_2''$ are piecewise smooth, we have
\[
w_0\in H^{1/2,2}(\mathbb{R}^+),\ \
w_1\in H^{3/2,2}(\mathbb{R}^+),\ \
w_2\in H^{5/2,2}(\mathbb{R}^+),
\]
where now the Nikolskii spaces on $\mathbb{R}^+$ are invoked, for example, see \cite{BL76,N75}. 
So the corresponding smoothness order of $w_k$ is fractional and equals $\lambda=k+1/2$, for $k=0,1,2$, and namely this is definitive for the convergence rates in the above chosen norms.

\par Consider a numerical method of the  $k$th approximation order for the IBVP for the wave equation. Below we compare the above 4th order semi-explicit compact scheme and the 2nd order classical explicit scheme 
\begin{gather*}
 \rho\Lambda_tv-L_hv=f\ \ \text{on}\ \ \omega_{\*h},
\\
 \rho\delta_tv^0-\tfrac{h_t}{2}L_hv^0=u_1+\tfrac{h_t}{2}f^0\ \ \text{on}\ \ \omega_{h},
\ \
v|_{\partial\omega_{\*h}}=g,\ \ v^0=u_0\ \ \text{on}\ \ \omega_{h}.
\end{gather*} 
Based on \cite{BL76,BTW75,Z94}, one can expect that, for $u_0$ of the Nikolskii smoothness order $\lambda$ and $u_1$ of the such smoothness order $\lambda-1$, the errors $e_{L_h^2}$ and  $e_{H_h^1}$ together with $e_{\mathcal{E}_h}$ have  the convergence rates, respectively,
\begin{equation}
p=\frac{k}{k+1}\lambda,\ \ \frac{k}{k+1}(\lambda-1),\ \ \text{for}\ \ 0\leq\lambda\leq k+1.
\label{theor orders}
\end{equation}
Moreover, $f=w_1$ is independent of $t$ and, thus, its 2nd order Sobolev derivatives $\partial_t^2f$ and $\partial_t\partial_rf$ are zero and $f=f|_{t=0}$.
Therefore, for such $f$ of the Nikolskii smoothness order $\lambda-2$ in $r$, one can expect the same convergence rates. 
Clearly both the rates increase as $k$ grows 
(and tend to $\lambda$ and $\lambda-1$, respectively) 
that ensures the advantage of higher-order methods over second-order ones not only for smooth data but for nonsmooth data as well, in contrast to the elliptic and parabolic cases.
Note that the difference between the two rates \eqref{theor orders} equals $k/(k+1)$ (i.e., 0.8 for $k=4$ and $0.666$... for $k=2$) independently on $\lambda$.

\par We recall how to construct the exact spherically symmetric solutions $u(x,t)=u(r,t)$ with $r=|x|$ for the Cauchy problem to the 3D wave equation with the spherically symmetric data.
They satisfy the IBVP on the half-line
\begin{gather*}
\partial_t^2u=a^2\frac{1}{r^2}\partial_r(r^2\partial_r u)+f(r,t)\ \ \text{for}\ \ r>0,\ t>0,
\\
\partial_ru|_{r=0}=0\ \ \text{for}\ \ t>0,\ \
u|_{t=0}=u_0(r),\ \ \partial_tu|_{t=0}=u_1(r)\ \ \text{for}\ \ r>0.
\end{gather*}
The function $v(r,t):=ru(r,t)$ solves the standard IBVP on the half-line
\begin{gather*}\partial_t^2v=a^2\partial_r^2v+rf(r,t)\ \ \text{for}\ \ r>0,\ t>0,
\\
v|_{r=0}=0\ \ \text{for}\ \ t>0,\ \ v|_{t=0}=ru_0(r),\ \ \partial_tv|_{t=0}=ru_1(r)\ \ \text{for}\ \ r>0. 
\end{gather*}
We extend $v$ oddly in $r$: 
$v(-r,t)=-v(r,t)$ for $r>0$, $t\geq 0$.
The extended function $v$ satisfies already the Cauchy problem for the 1D wave equation
\begin{gather*}
\partial_t^2v=a^2\partial_r^2v+rf(|r|,t)\ \ \text{for}\ \ r\in\mathbb{R},\ t>0,
\\
 v|_{t=0}=ru_0(|r|),\ \ \partial_tv|_{t=0}=ru_1(|r|)\ \ \text{for}\ \ r\in\mathbb{R}
\end{gather*}
with the oddly in $r$ extended data $u_0,u_1$ and $f$.
Thus, $u$ can be represented by the d'Alembert-type formula
\begin{gather}
u(r,t)=\frac{1}{2r}\big[(r-at)u_0(|r-at|)+(r+at)u_0(|r+at|)\big]
 +\frac{1}{2ar}\int_{r-at}^{r+at}q u_1(q)\,dq
\nonumber\\
+\frac{1}{2ar}\int_0^t
 \int_{r-a(t-\tau)}^{r+a(t-\tau)}q f(q,\tau)\,dq d\tau\ \ \text{for}\ \ r>0,\ t\geq 0.
 \label{exact sol}
\end{gather}

\par We deal with the following six cases: (a) $u_0=w_1$; (b) $u_0=w_2$; (c) $u_1=w_0$; (d) $u_1=w_1$; (e)  $f=w_0$; (f) 
$f=w_1$. Here we assume that the data unmentioned in the cases are zero, for example, $u_1=0$ and $f=0$ in case (a).
For brevity, we omit
the explicit lengthy formulas for $u(r,t)$ in these cases.

\par
\begin{figure}[tbh!]
\begin{minipage}{0.5\textwidth}
\center{\includegraphics[width=1\linewidth]{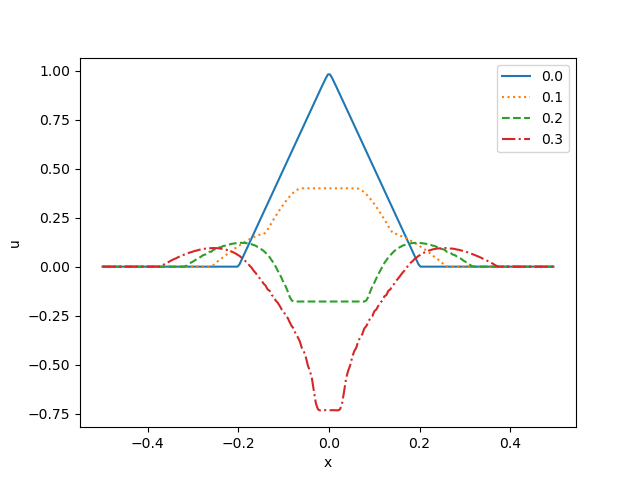}} (a) \\
 \end{minipage}
 \hfill
 \begin{minipage}{0.5\textwidth}
\center{\includegraphics[width=1\linewidth]{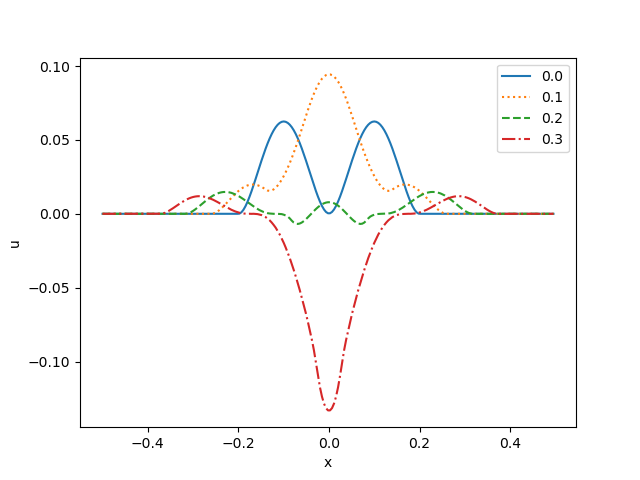}} (b) \\
 \end{minipage}
 \vfill
 \begin{minipage}{0.5\textwidth}
\center{\includegraphics[width=1\linewidth]{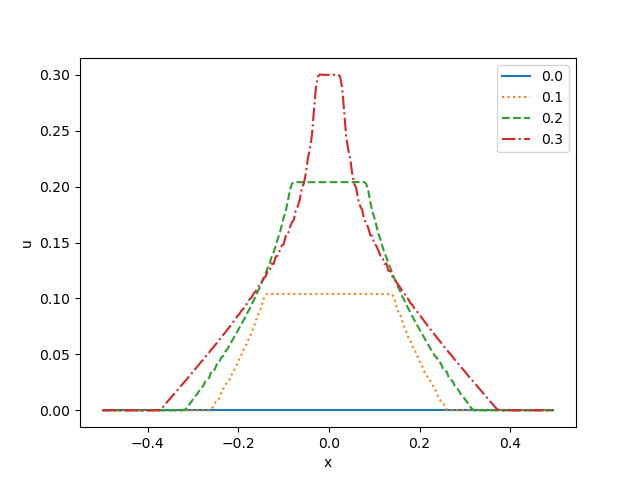}} (c) \\
 \end{minipage}
 \hfill
 \begin{minipage}{0.5\textwidth}
\center{\includegraphics[width=1\linewidth]{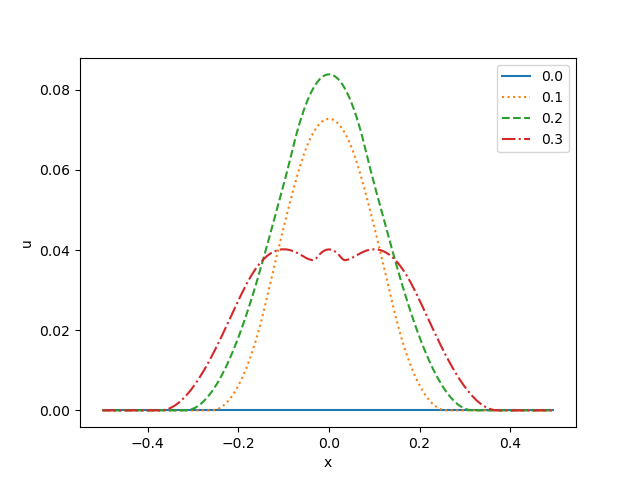}} (d) \\
 \end{minipage}
 \vfill
 \begin{minipage}{0.5\textwidth}
\center{\includegraphics[width=1\linewidth]{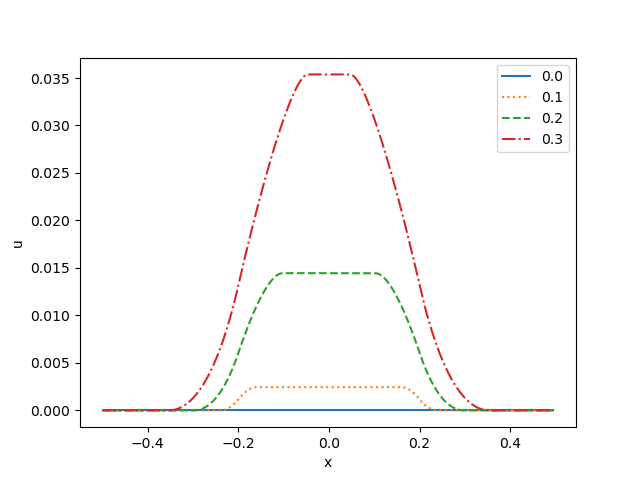}} (e) \\
 \end{minipage}
 \hfill
 \begin{minipage}{0.5\textwidth}
\center{\includegraphics[width=1\linewidth]{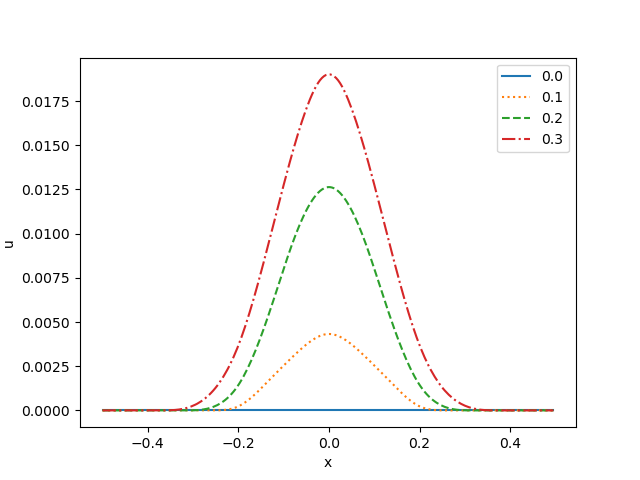}} (f) 
\\
 \end{minipage}
  \caption{Example 2. Numerical solutions in the section $y=z=0$ for $t=0,0.1,0.2$ and $0.3$ in the cases: (a) $u_0=w_1$; (b)~$u_0=w_2$; (c) $u_1=w_0$; (d) $u_1=w_1$; (e) $f=w_0$; (f)  $f=w_1$.}
\label{numer solut}
\end{figure}

 We take the symmetric cubic domain $\Omega=(-X/2,X/2)^3$ with $X=1$ and choose the parameters $a=1/\sqrt{3}$, $T=0.3$ and $r_0=0.2$. 
For such parameters, the exact solution is zero at the boundary for all $0\leq t\leq T$.
We take the same values of $N$ and $M$ as above for comparison. Since the values of $N$ are odd, we avoid the formal singularity at $r=0$ in summands of formula \eqref{exact sol} when computing $u$.
The numerical solutions in cases (a)-(f) are illustrated in Fig. \ref{numer solut} (visually they do not defer from the exact solutions).
As above, we compute the errors $e_{L_h^2}$,  $e_{H_h^1}$ and $e_{\mathcal{E}_h}$, the corresponding practical convergence rates
$p_{L_h^2}$, $p_{H_h^1}$ and $p_{\mathcal{E}_h}$ according to formula \eqref{pr_err_order} as well as the expected theoretical convergence rates according to formulas \eqref{theor orders} (recall that, in Tables, the last rates are formally attributed to the value $N=M=*$).

\par In Tables \ref{tab: u0=w1} and  \ref{tab: u0=w2}, we present them for the semi-explicit higher-order scheme ($k=4$) and classical explicit scheme ($k=2$) in the cases $u_0=w_1$ and $u_0=w_2$, respectively. 
Hereafter all the errors are many orders of magnitude worse than those for smooth data in Example 1. Also the error $e_{L_h^2}$ is much less (about two orders of magnitude) than $e_{H_h^1}$ and $e_{\mathcal{E}_h}$ compared with Example 1. 
In Table \ref{tab: u0=w1}, the errors $e_{H_h^1}$ and $e_{\mathcal{E}_h}$ are rather close to each other.
We see the very good agreement between $p_{L_h^2}$ and the corresponding theoretical convergence rates both for $k=4$ and $k=2$.
The values of $p_{H_h^1}$ and $p_{\mathcal{E}_h}$ are notably less than the theoretical rates, but at least for $k=4$ they are definitely higher than for $k=2$.
The latter is not surprising since, for low convergence rates, much higher values of $N$ and $M$ are often required to achieve better agreement, even in 1D case, cf. \cite{ZK18}. 
In Table \ref{tab: u0=w2}, the errors are about two orders of magnitude less than the corresponding ones in the previous Table, and now we observe the fine agreement between the practical and expected theoretical convergence rates.
\begin{table}[ht]
\begin{center}
\caption{Example 2(a). Errors $e_{L_h^2}$,  $e_{H_h^1}$ and $e_{\mathcal{E}_h}$ and convergence rates $p_{L_h^2}$, $p_{H_h^1}$ and $p_{\mathcal{E}_h}$, for $u_0=w_1$}
\label{tab: u0=w1}
\begin{tabular}{crrcccccc}
\hline\noalign{\smallskip}
& $N$ & $M$ & 
$e_{L_h^2}$ & $e_{H_h^1}$ & $e_{\mathcal{E}_h}$ &
$p_{L_h^2}$ & $p_{H_h^1}$ & $p_{\mathcal{E}_h}$ \\
\noalign{\smallskip}
\hline
\noalign{\smallskip}
$k=4$ & 81 & 27& 4.444064E--4 &8.412731E--2 &7.300600E--2 &$-$   &$-$    &$-$ \\
& 135& 45& 2.308279E--4 &7.058185E--2 &6.029682E--2 &1.282 &0.344 &0.374 \\
& 225& 75& 1.208428E--4 &5.859705E--2 &5.162431E--2 &1.267 &0.364 &0.304 \\
& 375&125& 6.468215E--5 &5.003797E--2 &4.327436E--2 &1.224 &0.309 &0.345 \\
& $*$&  $*$  &$-$   &$-$   &$-$           &1.200 &0.400 &0.400 \\
\hline
\noalign{\smallskip}
$k=2$ & 81 & 27& 9.032524E--4 &1.316278E--1 &1.104166E--1 &$-$   &$-$    &$-$ \\
      & 135& 45& 5.365153E--4 &1.168099E--1 &9.645128E--2 &1.020 &0.234 &0.265 \\
      & 225& 75& 3.176704E--4 &1.016564E--1 &8.528642E--2 &1.026 &0.272 &0.241 \\
      & 375&125& 1.899327E--4 &8.919430E--2 &7.411805E--2 &1.007 &0.256 &0.275 \\
& $*$&  $*$  &$-$    &$-$     &$-$              &1.000 &0.333 &0.333 \\
\hline
\end{tabular}
\end{center}
\end{table}
\begin{table}[ht]
\begin{center}
\caption{Example 2(b). Errors $e_{L_h^2}$,  $e_{H_h^1}$ and $e_{\mathcal{E}_h}$ and convergence rates $p_{L_h^2}$, $p_{H_h^1}$ and $p_{\mathcal{E}_h}$, for $u_0=w_2$}
\label{tab: u0=w2}
\begin{tabular}{crrcccccc}
\hline\noalign{\smallskip}
& $N$ & $M$ & 
$e_{L_h^2}$ & $e_{H_h^1}$ & $e_{\mathcal{E}_h}$ &
$p_{L_h^2}$ & $p_{H_h^1}$ & $p_{\mathcal{E}_h}$ \\
\noalign{\smallskip}
\hline
\noalign{\smallskip}
$k=4$ &81 & 27& 1.801824E--5 &3.162480E--3 &2.714343E--3 &$-$   &$-$    &$-$ \\
      &135& 45& 6.263696E--6 &1.740913E--3 &1.519010E--3 &2.068 &1.169 &1.136 \\
      &225& 75& 2.200167E--6 &9.725291E--4 &8.279443E--4 &2.048 &1.140 &1.188 \\
      &375&125& 7.757714E--7 &5.362174E--4 &4.555987E--4 &2.041 &1.165 &1.169 \\
      &$*$&  $*$  &$-$   &$-$   &$-$           &2.000 &1.200 &1.200 \\
\hline
\noalign{\smallskip}
$k=2$ &81 & 27& 7.455862E--5 &7.888876E--3 &6.577538E--3 &$-$   &$-$    &$-$ \\
      &135& 45& 3.209740E--5 &4.835505E--3 &4.021807E--3 &1.650 &0.9582 &0.9630 \\
      &225& 75& 1.379373E--5 &2.964015E--3 &2.440001E--3 &1.653 &0.9581 &0.9783 \\
      &375&125& 5.918562E--6 &1.804674E--3 &1.482635E--3 &1.656 &0.9713 &0.9752 \\
      &$*$&  $*$  &$-$    &$-$     &$-$        &1.667 &1.0000 &1.0000 \\
\hline
\end{tabular}
\end{center}
\end{table}

\par Tables \ref{tab: u1=w0} and \ref{tab: u1=w1} show the results for both the schemes in the cases $u_1=w_0$ and $u_1=w_1$, respectively. 
The agreement between the practical and expected theoretical convergence rates is good for
$p_{H_h^1}$ and $p_{\mathcal{E}_h}$ (despite that they are low in Table \ref{tab: u1=w0}) and not bad for $p_{L_h^2}$. In Table \ref{tab: u1=w1}, the agreement is better.
In addition, Fig.~\ref{fig:errors_time}(b) shows the dynamics of errors
$e_{L_h^2}$, $e_{H_h^1}$ and $e_{\mathcal{E}_h}$ in time in the second case;
the error $e_{L_h^2}$ is multiplied by 100 to make its behavior visible.
Once again we see a slow and rather close to linear growth of $e_{L_h^2}$. Now both $e_{H_h^1}$ and $e_{\mathcal{E}_h}$ grow much more rapidly but the growth is much slower than linear.
Notice oscillations in the behavior of all the errors; their amplitudes diminish in time.
Also the amplitudes diminish as $N$ and $M$ grows (in particular, for $N=81$ and $M=27$ they are larger; we do not show them). 
\begin{table}[ht]
\begin{center}
\caption{Example 2(c). Errors $e_{L_h^2}$, $e_{H_h^1}$ and $e_{\mathcal{E}_h}$ and convergence rates $p_{L_h^2}$, $p_{H_h^1}$ and $p_{\mathcal{E}_h}$, for $u_1=w_0$}
\label{tab: u1=w0}
\begin{tabular}{crrcccccc}
\hline\noalign{\smallskip}
      &$N$ & $M$ & 
$e_{L_h^2}$ & $e_{H_h^1}$ & $e_{\mathcal{E}_h}$ &
$p_{L_h^2}$ & $p_{H_h^1}$ & $p_{\mathcal{E}_h}$ \\
\noalign{\smallskip}
\hline
\noalign{\smallskip}
$k=4$ &81 & 27& 2.173270E--4 &3.947877E--2 &3.380764E--2 &$-$   &$-$    &$-$ \\
      &135& 45& 1.325706E--4 &3.207199E--2 &2.725321E--2 &0.968 &0.407 &0.422 \\
      &225& 75& 6.682470E--5 &2.591038E--2 &2.247714E--2 &1.341 &0.418 &0.377 \\
      &375&125& 3.317065E--5 &2.138626E--2 &1.834409E--2 &1.371 &0.376 &0.397 \\
      &$*$&  $*$  &$-$   &$-$   &$-$           &1.200 &0.400 &0.400 \\
\hline
\noalign{\smallskip}
$k=2$ &81 & 27& 3.542877E--4 &5.359390E--2 &4.470650E--2 &$-$   &$-$    &$-$ \\
      &135& 45& 2.165218E--4 &4.603917E--2 &3.790000E--2 &0.964 &0.297 &0.323 \\
      &225& 75& 1.233063E--4 &3.924349E--2 &3.274401E--2 &1.102 &0.313 &0.286 \\
      &375&125& 7.091202E--5 &3.370893E--2 &2.793154E--2 &1.083 &0.298 &0.311 \\
      &$*$&  $*$  &$-$    &$-$     &$-$        &1.000 &0.333 &0.333 \\
\hline
\end{tabular}
\end{center}
\end{table}
\begin{table}[ht]
\begin{center}
\caption{Example 2(d). Errors $e_{L_h^2}$,  $e_{H_h^1}$ and $e_{\mathcal{E}_h}$ and convergence rates $p_{L_h^2}$, $p_{H_h^1}$ and $p_{\mathcal{E}_h}$, for $u_1=w_1$}
\label{tab: u1=w1}
\begin{tabular}{crrcccccc}
\hline\noalign{\smallskip}
      &$N$ & $M$ & 
$e_{L_h^2}$ & $e_{H_h^1}$ & $e_{\mathcal{E}_h}$ &
$p_{L_h^2}$ & $p_{H_h^1}$ & $p_{\mathcal{E}_h}$ \\
\noalign{\smallskip}
\hline
\noalign{\smallskip}
$k=4$ &81 & 27& 3.486452E--6 &6.014188E--4 &5.177318E--4 &$-$   &$-$    &$-$ \\
      &135& 45& 1.201817E--6 &3.370573E--4 &2.921310E--4 &2.085 &1.134 &1.120 \\
      &225& 75& 4.077593E--7 &1.824089E--4 &1.552747E--4 &2.116 &1.202 &1.237 \\
      &375&125& 1.412191E--7 &9.888468E--5 &8.401369E--5 &2.076 &1.199 &1.202 \\
      &$*$&  $*$  &$-$   &$-$   &$-$           &2.000 &1.200 &1.200 \\
\hline
\noalign{\smallskip}
$k=2$ &81 & 27& 1.585246E--5 &1.435502E--3 &1.178713E--3 &$-$   &$-$    &$-$ \\
      &135& 45& 6.463676E--6 &8.736566E--4 &7.217399E--4 &1.756 &0.9721 &0.9602 \\
      &225& 75& 2.666898E--6 &5.271610E--4 &4.330681E--4 &1.733 &0.9890 &0.9999 \\
      &375&125& 1.110063E--6 &3.184605E--4 &2.613995E--4 &1.716 &0.9867 &0.9883 \\
      &$*$&  $*$  &$-$    &$-$     &$-$        &1.667 &1.0000 &1.0000 \\
\hline 
\end{tabular}
\end{center}
\end{table}

\par Tables \ref{tab: f=w_0} and \ref{tab: f=w_1} contain the results concerning the $H_h^1$- and $\mathcal{E}_h$-error norms for the both schemes in the case $f=w_0$ and $f=w_1$, respectively.
In Table \ref{tab: f=w_0}, all the errors are slightly smaller for $k=4$ than the corresponding ones for $k=2$. But the practical convergence rates do not differ so significantly and, for the least and other values of $N$ and $M$, they are closer to the theoretical ones for $k=4$ and $k=2$, respectively.  
In Table \ref{tab: f=w_1}, the agreement between practical and theoretical convergence rates is nice for $k=4$ but worse for $k=2$, but the practical rates for $k=4$ are definitely larger than for $k=2$.
In the both cases, the behavior of $e_{L_h^2}$ and $p_{L_h^2}$ seems still rather chaotic, and we omit them. 
\par Finally, in all the cases, we observe the advantages in error behavior of the semi-explicit higher order scheme over the standard explicit scheme.
\begin{table}[ht]
\begin{center}
\caption{Example 2(e). Errors $e_{H_h^1}$ and $e_{\mathcal{E}_h}$ and convergence rates $p_{H_h^1}$ and $p_{\mathcal{E}_h}$, for $f=w_0$}
\label{tab: f=w_0}
\vspace{1mm}
\begin{tabular}{crrcccc}
\hline\noalign{\smallskip}
      &$N$ & $M$ & $e_{H_h^1}$ & $e_{\mathcal{E}_h}$ &
$p_{H_h^1}$ & $p_{\mathcal{E}_h}$ \\
\noalign{\smallskip}
\hline
\noalign{\smallskip}
$k=4$ &81 & 27 &8.346416E--4 &5.933564E--4 &$-$    &$-$ \\
      &135& 45 &4.141449E--4 &3.080528E--4 &1.371 &1.283 \\
      &225& 75 &2.412233E--4 &1.795166E--4 &1.058 &1.057 \\
      &375&125 &1.451002E--4 &1.071878E--4 &0.995 &1.010 \\      &$*$&  $*$ &$-$  &$-$      &1.200 &1.200 \\
\hline
\noalign{\smallskip}
$k=2$ &81 & 27 &9.056331E--4 &6.841736E--4   &$-$    &$-$ \\
      &135& 45 &4.786699E--4 &3.766503E--4 &1.248 &1.168 \\
      &225& 75 &2.718633E--4 &2.150001E--4 &1.107 &1.098 \\
      &375&125 &1.586915E--4 &1.253125E--4 &1.054 &1.057 \\
      &$*$&  $*$ &$-$  &$-$      &1.000 &1.000 \\
\hline
\end{tabular}
\end{center}
\end{table}
\begin{table}[ht]
\begin{center}
\caption{Example 2(f). Errors $e_{H_h^1}$ and $e_{\mathcal{E}_h}$ and convergence rates $p_{H_h^1}$ and $p_{\mathcal{E}_h}$, for  $f=w_1$}
\label{tab: f=w_1}
\begin{tabular}{crrcccc}
\hline\noalign{\smallskip}
      &$N$ & $M$ & 
$e_{H_h^1}$ & $e_{\mathcal{E}_h}$ &
$p_{H_h^1}$ & $p_{\mathcal{E}_h}$ \\
\noalign{\smallskip}
\hline
\noalign{\smallskip}
$k=4$ &81 & 27 &5.459770E--6 &4.493218E--6   &$-$    &$-$ \\
      &135& 45 &1.857594E--6 &1.523632E--6 &2.111 &2.117 \\
      &225& 75 &6.319904E--7 &5.259791E--7 &2.111 &2.082 \\
      &375&125 &2.205933E--7 &1.838900E--7 &2.060 &2.057 \\
      &$*$&  $*$ &$-$  &$-$      &2.000 &2.000 \\
\hline
\noalign{\smallskip}
$k=2$ &81 & 27 &4.351183E--5 &3.003301E--5   &$-$    &$-$ \\
      &135& 45 &1.676865E--5 &1.173838E--5 &1.867 &1.839 \\
      &225& 75 &6.422439E--6 &4.604730E--6 &1.879 &1.832 \\
      &375&125 &2.507550E--6 &1.836635E--6 &1.841 &1.799 \\
      &$*$&  $*$ &$-$  &$-$      &1.666 &1.666 \\
\hline
\end{tabular}
\end{center}
\end{table}

\smallskip\par
\begin{figure}[tbh!]
\begin{minipage}{0.5\textwidth}
\center{\includegraphics[width=1\linewidth]{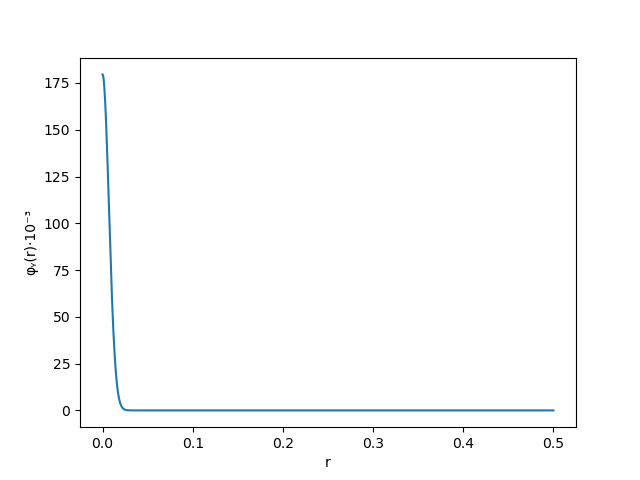}} (a) \\
 \end{minipage}
 \hfill
 \begin{minipage}{0.5\textwidth}
\center{\includegraphics[width=1\linewidth]{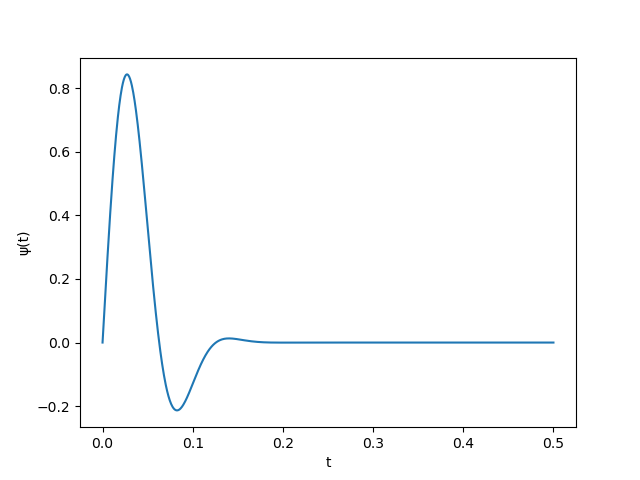}} (b) \\
\end{minipage}
  \caption{Example 3. Multipliers $\varphi_\gamma(r)$ for $\gamma=10000$ and $\psi(t)$ of the source function.}
\label{multipliers phi and psi} 
\end{figure}

{\bf Example 3.}
This example is similar to a part of Example 6 in \cite{ZC23}, where one computes the wave propagation in the three-layer medium in $x$, with
the different speeds of sound $1.5$, $1$ and $3$ km/sec, respectively, in its left, middle
and right layers of the domain $\Omega=(0,X)^3$ with $X=3$ km. 
The layers have  the same thickness 1 km. 
Recall that, in these layers, we set $a=1$ but take  $\rho=4/9, 1$ and $1/9$ (sec/km)$^2$, respectively. 
The source is the Ricker-type wavelet (popular in geophysics) smoothed in space
\[
f(x,y,z,t)=\varphi_\gamma(r)\psi(t),\ \ \text{with}\ \ \varphi_\gamma(r)=\Big(\frac{\pi}{\gamma}\Big)^{3/2}e^{-\gamma r^2},\ \ 
\int_{\mathbb{R}^3}\varphi_\gamma(r)\,dxdydz=1,
\ \ \psi(t)=\sin(50t)e^{-200t^2},
\]
where now $r=|(x-x_0,y-y_0,z-z_0)|$,  $(x_0,y_0,z_0)=(1.5,1.5,1.5)$ is the center of domain,  and $\gamma\gg 1$ is a parameter. Also we take $u_0=u_1=0$, 
$\gamma=10000$ and $T=0.8$ sec. 
The multipliers $\varphi_\gamma(r)$ and $\psi(t)$ are plotted in Fig. \ref{multipliers phi and psi}.

\par
\begin{figure}[tbh!]
\begin{minipage}{0.5\textwidth}
\center{\includegraphics[width=1\linewidth]{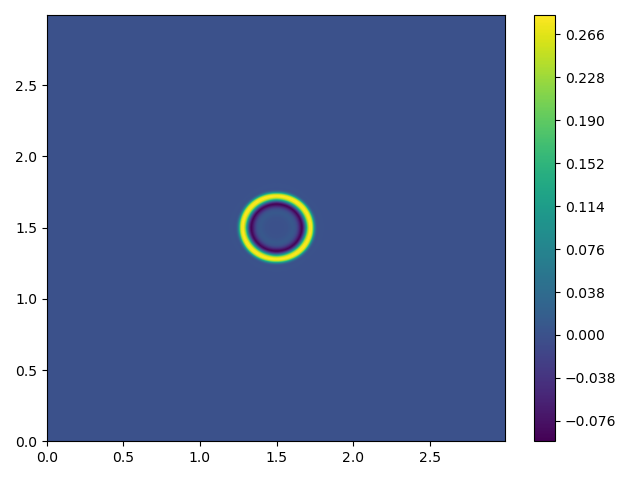}} (a) \\
\end{minipage}
 \hfill
\begin{minipage}{0.5\textwidth}
\center{\includegraphics[width=1\linewidth]{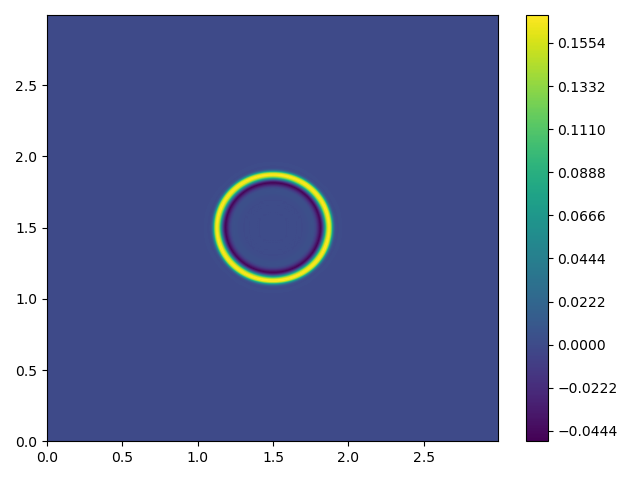}} (b) \\
\end{minipage}
 \vfill
\begin{minipage}{0.5\textwidth}
\center{\includegraphics[width=1\linewidth]{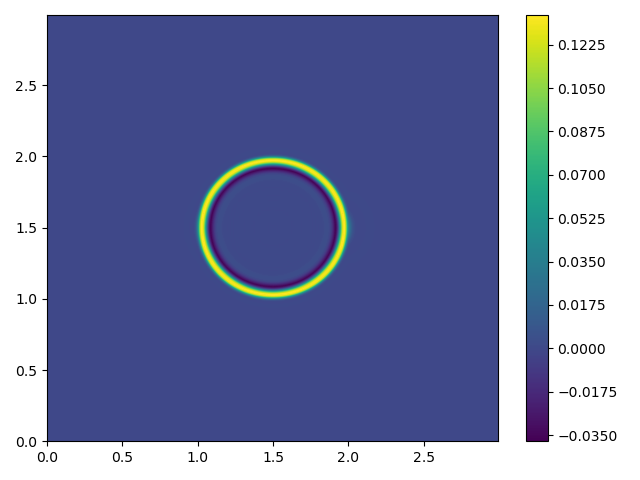}} (c) \\
\end{minipage}
 \hfill
\begin{minipage}{0.5\textwidth}
\center{\includegraphics[width=1\linewidth]{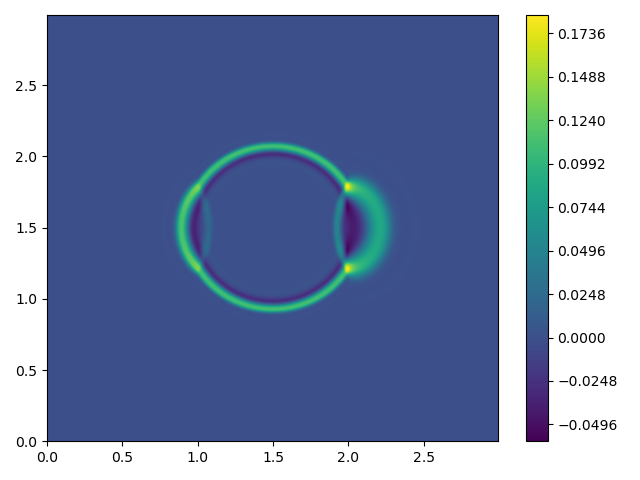}} (d) \\
\end{minipage}
 \vfill
\begin{minipage}{0.5\textwidth}
\center{\includegraphics[width=1\linewidth]{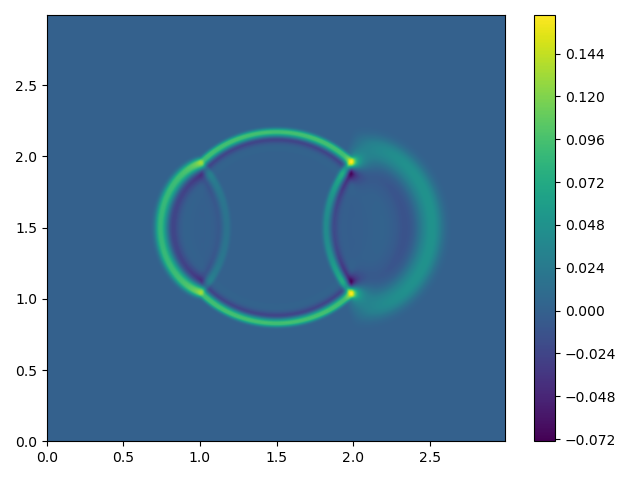}} (e) \\
\end{minipage}
 \hfill
\begin{minipage}{0.5\textwidth}
\center{\includegraphics[width=1\linewidth]{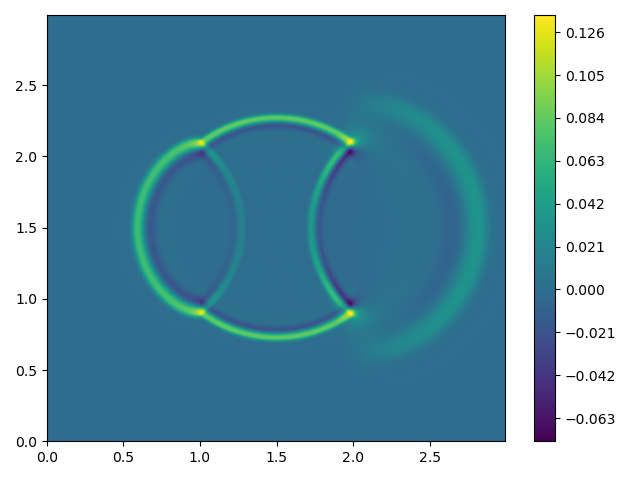}} (f) \\
\end{minipage}
\caption{Example 3. Contour levels of wavefields in section $z=1.5$ computed for $h=0.075$ and $h_t=1/700$ at: (a) $t=0.25$; (b) $t=0.4$; (c) $t=0.5$; (d) $t=0.6$; (e) $t=0.7$; (f) $t=0.8$.}
\label{contour figures}  
\end{figure}

 We present figures similar to those given in \cite{ZC23}. 
The results are in general close to those from \cite{ZC23} where $\varphi_\gamma(r)$ was taken as the delta-function concentrated at $(x_0,y_0,z_0)$, but with some essential differences that we discuss below.
In Fig. \ref{contour figures}, we demonstrate the 2D
wave fields in section $z=1.5$ computed for $N=400$ and $M=560$ (i.e., $h=0.0075$ and $h_t=1/700$) at six time moments.
In Fig. \ref{contour figures}(a) and (b), the spherical wave is expanding inside the middle layer and leaving the enlarging neighborhood of the center of the domain.
In Fig. \ref{contour figures}(c), the wave reaches the boundaries between the layers and only begins to penetrate into the right layer.
In Fig. \ref{contour figures}(d), the penetration to the left and right layers together with reflections from the boundaries of these layers inside the middle layer are already visible well. Both the effects are stronger for the right layer with the maximal speed of sound.
In Fig. \ref{contour figures}(e) and (f), both the penetrations and reflections evolve rapidly and, in the last figure, the wave front almost approaches the right boundary of the domain. 

\par
\begin{figure}[tbh!]
\begin{minipage}{0.5\textwidth}
\center{\includegraphics[width=1\linewidth]{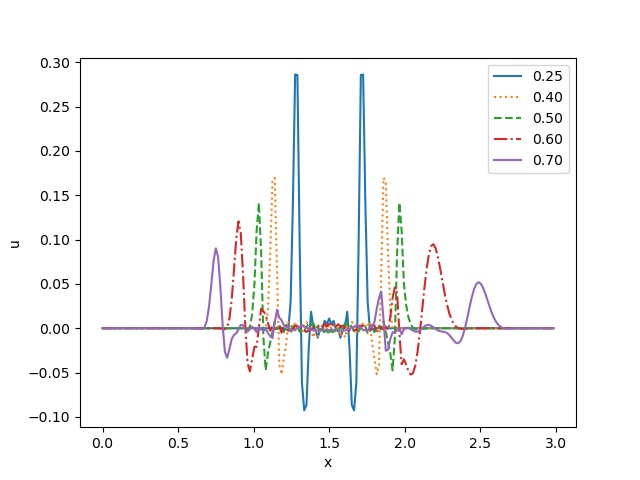}} (a) \\
 \end{minipage}
 \hfill
 \begin{minipage}{0.5\textwidth}
\center{\includegraphics[width=1\linewidth]{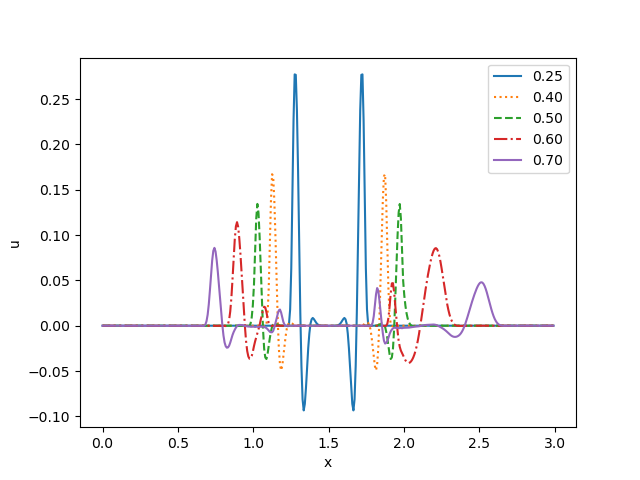}} (b) \\
\end{minipage}
  \caption{Example 3. Dynamics of the waves at different times for $y=z=1.5$: (a) $h=0.015$, $h_t=1/352$ for $t=0.25$ and $h_t=1/350$ for other times; (b) $h=0.075$ and $h_t=1/700$.}
\label{dyn at fixed y z} 
\end{figure}

\begin{figure}[tbh!]
\begin{minipage}{0.5\textwidth}
\center{\includegraphics[width=1\linewidth]{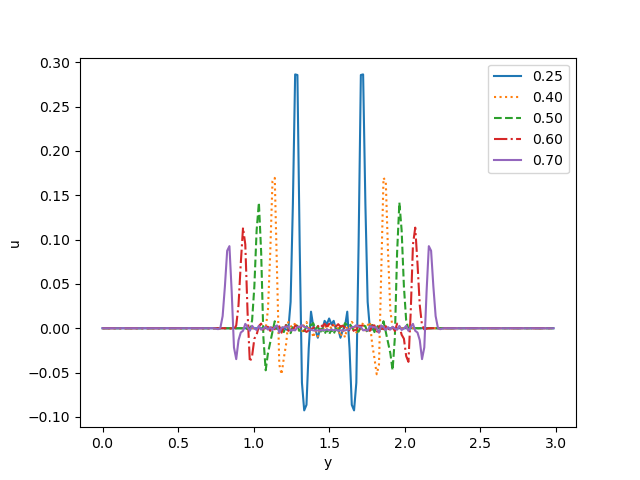}} (a) \\
 \end{minipage}
 \hfill
 \begin{minipage}{0.5\textwidth}
\center{\includegraphics[width=1\linewidth]{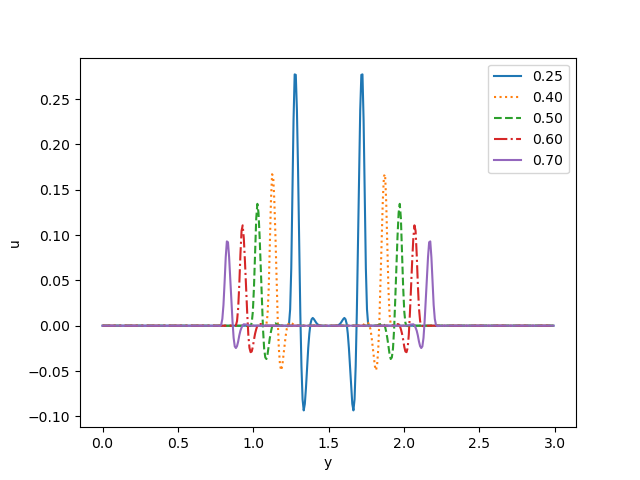}} (b) \\
\end{minipage}
\caption{Example 3. Dynamics of the waves at different times for $x=z=1.5$:
  (a) $h=0.015$, $h_t=1/352$ for $t=0.25$ and $h_t=1/350$ for other times; (b) $h=0.0075$ and $h_t=1/700$.}
\label{dyn at fixed x z}
 \end{figure}

 We supplement the last figure by Fig. \ref{dyn at fixed y z} and \ref{dyn at fixed x z} showing the 1D dynamics of the wave at different times for $y=z=1.5$
and $x=z=1.5$, respectively, at five (from above six) time moments. These figures are nonsymmetric in $x$ and symmetric in $y$, respectively, as it should be.
Of course, this dynamics corresponds to those shown in the previous figure in the perpendicular directions $x=1.5$ and $z=1.5$, respectively, but it also contains additional information about the amplitude and form of the moving wave fronts. 
The results are given for $N=200$ as in \cite{ZC23} and $N=400$ as above, respectively.
In general, they are close, but the results for $N=200$ suffer from parasitic numerical oscillations. Unfortunately, the results in \cite{ZC23} (including the counterparts of Fig. \ref{contour figures}(e) and (f)) also contain such parasitic oscillations essentially distorting the exact solution in some subdomains where it is rather small.

\par We also make a note on the quality of contour levels.
Fig. \ref{contour figures} is plotted via the \textit{matplotlib.pyplot} library for Python 3.10.9. 
It is crucial to choose properly its parameter \textit{levels} defining the maximal amount of levels in use.
For $N=400$ and small values of \textit{levels} like 14, we encountered false visible oscillations and the false observation of the wave achieving the right boundary of the domain, without its reflection.
Setting that parameter like $N$ provided us with correct results.
 
\section*{Acknoledgemnents}
The work was supported by the Russian Science Foundation, grant no. 23-21-00061.

\section*{Declaration of competing interest}
The authors have no conflicts of interest to declare that are relevant to the content of this paper.

\end{document}